\documentclass[a4paper,11pt,thmsb]{article}
\usepackage{graphicx}
\usepackage{amsmath}
\usepackage{amsfonts}
\usepackage{amssymb}
\newtheorem{theorem}{Th\'eor\`eme}

\newtheorem{corollary}{Corollaire}
\newtheorem{definition}{D\'efinition}

\newtheorem{lemma}{Lemme}

\newenvironment{proof}[1][Proof]{\textbf{#1.} }{\ \rule{0.5em}{0.5em}}
\begin{document}

\title{Hermitian forms, trace equations\\and application to codes}
\date{November 14$^{th}$, 2001}
\author{Dany-Jack Mercier\bigskip\thanks{NT/0111191 [ccof0007] v 1.00 Equipe
Applications de l'Alg\`{e}bre et de l'Arithm\'{e}tique, D\'{e}partement de
Math\'{e}matiques \& Informatique, Universit\'{e} des Antilles-Guyane, Campus
Fouillole, F97159 Pointe-\`{a}-Pitre cedex, France, e-mail~:~dany-jack.mercier@univ-ag.fr.}}
\maketitle

%

\parindent0cm%

\begin{center}
\textbf{Abstract:}
\end{center}

\begin{quote}
We provide a systematic study of sesquilinear hermitian forms and a new proof
of the calculus of some exponential sums defined with quadratic hermitian
forms. The computation of the number of solutions of equations such as
$\operatorname*{Tr}\nolimits_{\mathbb{F}_{t}/\mathbb{F}_{s}}(f\left(
x\right)  +v.x)=0$ or $\operatorname*{Tr}\nolimits_{\mathbb{F}_{t}%
/\mathbb{F}_{s}}(f\left(  x\right)  )=a$ allows us to construct codes and to
obtain their parameters.
\end{quote}

\begin{center}
\textbf{Key-Words:}
\end{center}

\begin{quote}
Sesquilinear Hermitian Forms, Quadratic Forms, Finite Fields, Traces,
Exponential Sums, Codes.
\end{quote}

\begin{center}
\textbf{MSC-Class:}
\end{center}

\begin{quote}
94B05, 11T23, 11T71.
\end{quote}

\begin{center}
\textbf{Versions:}
\end{center}

\begin{quote}
Minor changes in version 2 (21$^{th}$ april 2003).
\end{quote}

\section{Introduction}

Let $\mathbb{F}_{t}$ be the finite field with $t$ elements and characteristic
$p$. Our first purpose is to adapt the classical hermitian form Theory on
$\mathbb{C}$ to the case of the finite field $\mathbb{F}_{t^{2}}$, considering
the involution $x\mapsto x^{t}$ in $\mathbb{F}_{t^{2}}$ instead of the
application $z\mapsto\overline{z}$. Then we introduce exponential sums
associated with quadratic hermitian forms and obtain the number of solutions
of certain equations on $\mathbb{F}_{t}$. At this point it is easy to
construct two linear codes using the same method as Reed-Muller codes, to get
their parameters and to compare them to the classical Reed-Muller construction.\smallskip

A general introduction to hermitian forms over a finite field is given by Bose
and Chakravarti in~\cite{bose}, and the use of those objects in coding Theory
has been discussed for instance in \cite{cherdieu}, \cite{mercier 1997}
or~\cite{mercier 2001}.\smallskip

My first contribution will be to review in details the main results
of~\cite{bose} and~\cite{cherdieu}, giving alternative proofs of some
important results in the Theory. It is worth pointing out that the existence
of a $H$-orthogonal basis is shown by induction without explicit calculus on
rows and columns of a matrix as in~\cite{bose} (see Theorem~\ref{990922a}).
Another example is given by a new and staightforward proof of the most
precious result of the paper of Cherdieu~\cite{cherdieu} (see
Theorem~\ref{cher}).\smallskip

Section~\ref{991014a} is devoted to the construction of the code $\Gamma$ of
\cite{cherdieu}, and Section~\ref{991101a} use the same argument to give
another example of such a linear code. As this paper wants to remain
self-contained, an annex in Section~\ref{991013b} will summarize without
proofs the relevant material on characters on a finite group.

\section{Sesquilinear hermitian forms on $\mathbb{F}_{t^{2}}^{N}$}

Let $N$ be an integer $\geq1$ and let $E$ be the vector space $\mathbb{F}%
_{t^{2}}^{N}$.

\begin{definition}
A function $H:\mathbb{F}_{t^{2}}^{N}\times\mathbb{F}_{t^{2}}^{N}$
$\rightarrow\mathbb{F}_{t^{2}}$ is a \textbf{sesquilinear form on}%
\emph{\ }$E=\mathbb{F}_{t^{2}}^{N}$ if it is semi-linear in the first variable
and linear in the second variable, i.e.
\[%
\begin{array}
[c]{l}%
\left(  1\right)  \quad\forall\lambda,\mu\in\mathbb{F}_{t^{2}}\quad\forall
x,x^{\prime},y\in E\quad H\left(  \lambda x+\mu x^{\prime},y\right)
=\lambda^{t}H\left(  x,y\right)  +\mu^{t}H\left(  x^{\prime},y\right)  ,\\
\left(  2\right)  \quad\forall\lambda,\mu\in\mathbb{F}_{t^{2}}\quad\forall
x,y,y^{\prime}\in E\quad H\left(  x,\lambda y+\mu y^{\prime}\right)  =\lambda
H\left(  x,y\right)  +\mu H\left(  x,y^{\prime}\right)  .
\end{array}
\]
The sesquilinear form $H$ is called\textbf{\ hermitian} if
\[
\left(  3\right)  \quad\forall x,y\in E\quad H\left(  x,y\right)  =H\left(
y,x\right)  ^{t}.
\]
A sesquilinear hermitian form will be called a \textbf{hermitian form on }$E$.
The vector space of all hermitian forms on $\mathbb{F}_{t^{2}}^{N}$ will be
denoted by$~\operatorname*{H}\left(  \mathbb{F}_{t^{2}}^{N}\right)  $.
\end{definition}

Note that properties $\left(  2\right)  $ and $\left(  3\right)  $ give
$\left(  1\right)  $, and that $H\left(  x,x\right)  \in\mathbb{F}_{t}$ for
all $x\in E$ as soon as $H$ is a hermitian~form.\smallskip

If $x\in\mathbb{F}_{t^{2}}$ we put $\overline{x}=x^{t}$ and we say that
$\overline{x}$ is the conjugate of$~x$. If $\alpha\in\mathbb{F}_{t^{2}}$
satisfies $\mathbb{F}_{t^{2}}=\mathbb{F}_{t}\left(  \alpha\right)  $, each
element $x$ of $\mathbb{F}_{t^{2}}$ is uniquely written as $x=a+b\alpha$ with
$a$ and $b$ in~$\mathbb{F}_{t}$. Then $\overline{x}=\left(  a+b\alpha\right)
^{t}=a+b\alpha^{t}=a+b\overline{\alpha}.$

\begin{definition}
Let $A=\left(  a_{ij}\right)  _{i,j}$ be a square matrix with $i,j=1,...,N$
and with entries in$~\mathbb{F}_{t^{2}}$. We denote by $\overline{A}$ the
matrix $\overline{A}=\left(  \overline{a}_{ij}\right)  _{i,j}$ and by $^{T}A$
the transpose of $A$. The \textbf{conjugate} of $A=\left(  a_{ij}\right)
_{i,j}$ is the matrix $A^{\ast}=\,^{T}\left(  \overline{A}\right)  =\left(
a_{ji}^{t}\right)  _{i,j}$. The matrix $A$ is \textbf{hermitian} if~$A^{\ast
}=A$.
\end{definition}

If $e=\left(  e_{1},...,e_{N}\right)  $ is a basis of $E$ and if $H$ is a
sesquilinear form on $E$,
\[
H\left(  x,y\right)  =H\left(  \sum\limits_{i}x_{i}e_{i},\sum\limits_{j}%
y_{j}e_{j}\right)  =\sum\limits_{i,j}\overline{x}_{i}y_{j}H\left(  e_{i}%
,e_{j}\right)  =X^{\ast}MY
\]
where $M=\left(  H\left(  e_{i},e_{j}\right)  \right)  _{i,j}$, $X=\,^{T}%
\left(  x_{1},...,x_{N}\right)  $ and $Y=\,^{T}\left(  y_{1},...,y_{N}\right)
$. We say that $M$ is \textbf{the matrix of }$H$\textbf{\ in the basis~}$e$,
and we write$~M=\operatorname*{Mat}\left(  H;e\right)  $.

\begin{theorem}
A sesquilinear form $H$ is hermitian if, and only if, its matrix
$\operatorname*{Mat}\left(  H;e\right)  $ in a basis $e$ is hermitian.
\end{theorem}

\textbf{Proof }: Let $M$ denotes the matrix $\operatorname*{Mat}\left(
H;e\right)  $. If $H$ is a hermitian form, $H\left(  e_{i},e_{j}\right)
=\overline{H\left(  e_{j},e_{i}\right)  }$ implies $M^{\ast}=M$. Conversely,
$M^{\ast}=M$ implies
\[
\forall X,Y\quad\overline{H\left(  Y,X\right)  }=\overline{Y}^{\ast}%
\overline{M}\,\overline{X}=\,^{T}Y\,\overline{M}\,\overline{X}=\,^{T}\left(
^{T}Y\,\overline{M}\,\overline{X}\right)  =X^{\ast}MY=H\left(  X,Y\right)
\text{.%
\endproof
}%
\]

\begin{corollary}
$\dim_{\mathbb{F}_{t}}\operatorname*{H}\left(  \mathbb{F}_{t^{2}}^{N}\right)
=N^{2} $.
\end{corollary}

\textbf{Proof }: The matrix $H=\left(  a_{ij}\right)  $ of a hermitian form
depends on $\frac{N^{2}-N}{2}$ coefficients $a_{ij}$ (where $1\leq j<i\leq N$)
in $\mathbb{F}_{t^{2}}$ and $N$ coefficients $a_{ii}$ ($1\leq i\leq N$)
in~$\mathbb{F}_{t}$. Thus we have $2\times\frac{N^{2}-N}{2}+N=N^{2}$
independent parameters.%
\endproof
\medskip

Let $P_{e}^{e^{\prime}}=P$ denotes a change of coordinates from a basis
$e=\left(  e_{1},...,e_{N}\right)  $ to another basis $e^{\prime}=\left(
e_{1}^{\prime},...,e_{N}^{\prime}\right)  $. Let $X=\,^{T}\left(
x_{1},...,x_{N}\right)  $ and $X^{\prime}=\,^{T}\left(  x_{1}^{\prime
},...,x_{N}^{\prime}\right)  $ stands for the $N$-tuplets of coordinates of
the same vector in basis $e$ and $e^{\prime}$. Then
\[
H\left(  X,Y\right)  =X^{\ast}MY=\left(  PX^{\prime}\right)  ^{\ast}M\left(
PY^{\prime}\right)  =X^{\prime\ast}\left(  P^{\ast}MP\right)  Y^{\prime}%
\]
and $\operatorname*{Mat}\left(  H;e^{\prime}\right)  =P^{\ast}MP$ is the
matrix of $H$ in the new basis~$e^{\prime}$.

\section{Kernel and rank of $H$}

If $x\in E$ and if $H$ denotes a hermitian form, we define the linear
application $H\left(  x,.\right)  $ in the dual $E^{\ast}$ by:
\[%
\begin{array}
[c]{cccc}%
H\left(  x,.\right)  : & E & \longrightarrow & \mathbb{F}_{t^{2}}\\
& y & \longmapsto &  H\left(  x,y\right)  .
\end{array}
\]
The map:
\[%
\begin{array}
[c]{cccc}%
\widetilde{H}: & E & \longrightarrow &  E^{\ast}\\
& x & \longmapsto &  H\left(  x,.\right)
\end{array}
\]
is semi-linear, i.e. satisfies $\widetilde{H}\left(  \lambda x+x^{\prime
}\right)  =\overline{\lambda}\widetilde{H}\left(  x\right)  +\widetilde
{H}\left(  x^{\prime}\right)  $ for all vectors $x$, $x^{\prime}$ and all
$\lambda\in\mathbb{F}_{t^{2}}$. The extern law $\bullet$ defined by
$\lambda\bullet l=\overline{\lambda}.l$ gives us a new vector space structure
on $E^{\ast}$. For convenience, we shall write $\overline{E^{\ast}}$ instead
of $E^{\ast}$ when we use this new extern law. The map $\widetilde
{H}:E\rightarrow E^{\ast}$ is semi-linear if, and only if, $\widetilde
{H}:E\rightarrow\overline{E^{\ast}}$ is linear, and we can introduce the
matrix of $\widetilde{H}$ in the basis $e=\left(  e_{1},...,e_{N}\right)  $ in
$E$ and the dual basis $e^{\ast}=\left(  e_{1}^{\ast},...,e_{N}^{\ast}\right)
$ in~$\overline{E^{\ast}}$. The linearity of $\widetilde{H}$ gives us the same
results as in the case of symmetric bilinear forms. Namely:

\begin{theorem}
The equality $\operatorname*{Mat}\left(  \widetilde{H};e,e^{\ast}\right)
=\operatorname*{Mat}\left(  H;e\right)  $ holds for all hermitian form$~H$.
\end{theorem}

\textbf{Proof }: Assume that $\operatorname*{Mat}\left(  \widetilde
{H};e,e^{\ast}\right)  =\left(  a_{ij}\right)  $. Then $\widetilde{H}\left(
e_{j}\right)  =\sum\limits_{i}a_{ij}\bullet e_{i}^{\ast}$ and
\[
\widetilde{H}\left(  e_{j}\right)  \left(  e_{k}\right)  =H\left(  e_{j}%
,e_{k}\right)  =\overline{a}_{kj}.
\]
Thus $\operatorname*{Mat}\left(  \widetilde{H};e,e^{\ast}\right)
=\,^{T}\,\overline{\operatorname*{Mat}\left(  H;e\right)  }%
=\operatorname*{Mat}\left(  H;e\right)  $.%
\endproof

\begin{definition}
The \textbf{kernel} $\operatorname*{Ker}H$ (resp. \textbf{rank}
$\operatorname*{rk}H$) of $H$ is the \textbf{kernel} (resp. \textbf{rank})
of~$\widetilde{H}$. Thus $\operatorname*{Ker}H=\left\{  x\in E\,/\,\forall
y\in E\quad H\left(  x,y\right)  =0\right\}  \;$and $\operatorname*{rk}%
H=\operatorname*{rk}\left(  \operatorname*{Mat}\left(  H;e\right)  \right)  $.
\end{definition}

\section{Orthogonality}

\begin{definition}
Let $H$ denotes a hermitian form on $E$. Vectors $x$ and $y$ are
\textbf{orthogonal} if $H\left(  x,y\right)  =0$. If $F$ is a subset of$~E$,
the subspace $F^{\perp}=\left\{  x\in E\,/\,\forall y\in F\quad H\left(
x,y\right)  =0\right\}  $ is called \textbf{the orthogonal of~}$F$.
\end{definition}

It is easily seen that:

\begin{theorem}
\label{orth}For all subspaces $F$ and $G$ in $E$,
\[%
\begin{array}
[c]{ll}%
F\subset\left(  F^{\bot}\right)  ^{\bot}, & \quad F\subset G\Rightarrow
G^{\bot}\subset F^{\bot},\\
\left(  F+G\right)  ^{\bot}=F^{\bot}\cap G^{\bot}, & \quad\left(  F\cap
G\right)  ^{\bot}\supset F^{\bot}+G^{\bot}.
\end{array}
\]
\end{theorem}

\begin{definition}
A basis $e=\left(  e_{1},...,e_{N}\right)  $ is \textbf{orthogonal} for the
hermitian form $H$ (we say $H$\textbf{-orthogonal}) if $H\left(  e_{i}%
,e_{j}\right)  =0$ when $i\neq j$. This means that the matrix
$\operatorname*{Mat}\left(  H;e\right)  $ is diagonal.
\end{definition}

\begin{lemma}
\label{990921a}Let $H$ denotes a sesquilinear form on $E$. If $t$ is odd and
if $q\left(  x\right)  =H\left(  x,x\right)  $, then for all $x,y$ in
$E$,\newline 1) $H\left(  x,y\right)  +H\left(  y,x\right)  =\frac{1}%
{2}\left[  q\left(  x+y\right)  -q\left(  x-y\right)  \right]  ,$\newline 2)
$H\left(  x,y\right)  -H\left(  y,x\right)  =\frac{1}{\alpha-\overline{\alpha
}}\left[  q\left(  x+\alpha y\right)  -q\left(  x+\overline{\alpha}y\right)
\right]  ,$\newline 3) $H\left(  x,y\right)  =\frac{1}{4}\left[  q\left(
x+y\right)  -q\left(  x-y\right)  \right]  +\frac{1}{2\left(  \alpha
-\overline{\alpha}\right)  }\left[  q\left(  x+\alpha y\right)  -q\left(
x+\overline{\alpha}y\right)  \right]  .$
\end{lemma}

\textbf{Proof }: We have
\begin{align*}
q\left(  x+y\right)  -q\left(  x-y\right)   &  =q\left(  x\right)  +H\left(
x,y\right)  +H\left(  y,x\right)  +q\left(  y\right) \\
&  \quad\quad\quad\quad\quad\quad-\left[  q\left(  x\right)  +H\left(
x,-y\right)  +H\left(  -y,x\right)  +H\left(  -y,-y\right)  \right] \\
&  =2\left(  H\left(  x,y\right)  +H\left(  y,x\right)  \right)
\end{align*}
and
\begin{align*}
q\left(  x+\alpha y\right)  -q\left(  x+\overline{\alpha}y\right)   &
=q\left(  x\right)  +\alpha H\left(  x,y\right)  +\alpha^{t}H\left(
y,x\right)  +\alpha^{t+1}q\left(  y\right) \\
&  \quad\quad\quad-\left[  q\left(  x\right)  +\alpha^{t}H\left(  x,y\right)
+\alpha H\left(  y,x\right)  +\alpha^{t+1}q\left(  y\right)  \right] \\
&  =\left(  \alpha-\alpha^{t}\right)  \left(  H\left(  x,y\right)  -H\left(
y,x\right)  \right)  .
\end{align*}
As $\alpha\notin\mathbb{F}_{t}$, we have $\alpha^{t}\neq\alpha$ and the Lemma
follows $%
\endproof
$

\begin{lemma}
\label{990921b}The norm application
\[%
\begin{array}
[c]{cccc}%
N_{\mathbb{F}_{t^{m}}/\mathbb{F}_{t}}: & \mathbb{F}_{t^{m}}^{\ast} &
\rightarrow & \mathbb{F}_{t}^{\ast}\\
& x & \mapsto &  x^{t^{m-1}+...+t+1}%
\end{array}
\]
is a multiplicative group epimorphism, and $\left|  N_{\mathbb{F}_{t^{m}%
}/\mathbb{F}_{t}}^{-1}\left(  b\right)  \right|  =\frac{t^{m}-1}{t-1}$ for
all$~b\in\mathbb{F}_{t}^{\ast}$.
\end{lemma}

\textbf{Proof }: It is obvious that the map $N_{\mathbb{F}_{t^{m}}%
/\mathbb{F}_{t}}$ is a morphism. Consider a primitive element $\alpha$
in$~\mathbb{F}_{t^{m}}$. It means that $\alpha$ generates the multiplicative
group$~\mathbb{F}_{t^{m}}^{\ast}$, hence
\[
\mathbb{F}_{t^{m}}=\left\{  0,1,\alpha,\alpha^{2},...,\alpha^{t^{m}%
-2}\right\}  .
\]
Each element $x^{t^{m-1}+...+t+1}$ lies in $\mathbb{F}_{t}$ as $\left(
x^{t^{m-1}+...+t+1}\right)  ^{t}=x^{t^{m-1}+...t+1}$, and for $1\leq u\leq
t-1$, all elements $\alpha^{u\left(  t^{m-1}+...+t+1\right)  }$ are different.
Thus
\[
\mathbb{F}_{t}=\left\{  0,\alpha^{\left(  t^{m-1}+...+t+1\right)  }%
,\alpha^{2\left(  t^{m-1}+...+t+1\right)  },...,\alpha^{\left(  t-1\right)
\left(  t^{m-1}+...+t+1\right)  }\right\}
\]
and $N_{\mathbb{F}_{t^{m}}/\mathbb{F}_{t}}$ is surjective. The decomposition
of the morphism $N_{\mathbb{F}_{t^{m}}/\mathbb{F}_{t}}$ gives
\[
\mathbb{F}_{t^{m}}^{\ast}/\operatorname*{Ker}\left(  N_{\mathbb{F}_{t^{m}%
}/\mathbb{F}_{t}}\right)  \simeq\mathbb{F}_{t}^{\ast},
\]
hence $\left|  \operatorname*{Ker}\left(  N_{\mathbb{F}_{t^{m}}/\mathbb{F}%
_{t}}\right)  \right|  =\frac{t^{m}-1}{t-1}$. If $a\in N_{\mathbb{F}_{t^{m}%
}/\mathbb{F}_{t}}^{-1}\left(  b\right)  $, then
\[
N_{\mathbb{F}_{t^{m}}/\mathbb{F}_{t}}\left(  x\right)  =b\Leftrightarrow
N_{\mathbb{F}_{t^{m}}/\mathbb{F}_{t}}\left(  xa^{-1}\right)  =1\Leftrightarrow
x\in a\operatorname*{Ker}\left(  N_{\mathbb{F}_{t^{m}}/\mathbb{F}_{t}%
}\right)
\]
and we deduce $\left|  N_{\mathbb{F}_{t^{m}}/\mathbb{F}_{t}}^{-1}\left(
b\right)  \right|  =\left|  \operatorname*{Ker}N_{\mathbb{F}_{t^{m}%
}/\mathbb{F}_{t}}\right|  =\frac{t^{m}-1}{t-1}$.%
\endproof

\begin{theorem}
\label{990922a}\textbf{Existence of }$H$-\textbf{orthogonal basis}\newline If
$t$ is odd and if $H$ is a hermitian form on$~E=\mathbb{F}_{t^{2}}^{N}%
$,\newline 1) We can find a $H$-orthogonal basis $e=\left(  e_{1}%
,...,e_{N}\right)  $, and assume $H\left(  e_{i},e_{i}\right)  =0$ or $1$ for
all$~i$,\newline 2) The number $r$ of non zero entries in the diagonal of
$\operatorname*{Mat}\left(  H;e\right)  $ is an invariant that depends only
on$~H$. It is the rank of$~H$.
\end{theorem}

\textbf{Proof }: The proof of 1) is by induction on ~$N$. Assume $N=1$. The
result is obvious if $H\left(  x,x\right)  =0$ for all$~x$. If $x\in E$
satisfies $H\left(  x,x\right)  =b\in\mathbb{F}_{t}^{\ast}$,
Lemma~\ref{990921b} shows the existence of $a\in\mathbb{F}_{t^{2}}^{\ast}$
such that~$a^{t+1}=b$. Hence $H\left(  \frac{x}{a},\frac{x}{a}\right)  =1$ and
we take the basis$~e_{1}=\frac{x}{a}$.\smallskip

Assuming the result holds for $N$, we will prove it for $\dim E=N+1$. We need
only consider two cases~:

\quad\quad- If $H\left(  x,x\right)  =0$ for all $x$, then all basis are
$H$-orthogonal by formula~3) in Lemma~\ref{990921a}.

\quad\quad- If there exists $x$ with $H\left(  x,x\right)  \neq0$, we proceed
as in the case $N=1$ to show the existence of $a\in\mathbb{F}_{t^{2}}^{\ast}$
such that $H\left(  \frac{x}{a},\frac{x}{a}\right)  =1$. The subspace
\[
F=\left(  Kx\right)  ^{\bot}=\left\{  y\in E\,/\,H\left(  x,y\right)
=0\right\}  =\operatorname*{Ker}\widetilde{H}\left(  x\right)
\]
is an hyperplane of $E$ as it is the kernel of a non zero linear form. But
$E=F\oplus\operatorname*{Vect}\left(  x\right)  $ since $x\notin F$, and the
induction hypothesis gives a $H$-orthogonal basis $\left(  e_{1}%
,...,e_{N-1}\right)  $ of $F$ with $H\left(  e_{i},e_{i}\right)  =0$ or $1$
for all$~i$. We check at once that $\left(  e_{1},...,e_{N-1},\frac{x}%
{a}\right)  $ is a $H$-orthogonal basis of $E$ and this complete the proof of
1). The second part of the Theorem follows from%
\[
\operatorname*{rk}H=\operatorname*{rk}\widetilde{H}=\operatorname*{rk}%
\operatorname*{Mat}\left(  H;e\right)  .%
\endproof
\]

\begin{definition}
A hermitian form $H$ is \textbf{non degenerate} if $E^{\perp}=\left\{
0\right\}  $.
\end{definition}

\begin{theorem}
Let $H$ be a hermitian form. The following conditions are equivalent:\newline
1) $H$ is non degenerate,\newline 2) $\operatorname*{Ker}\widetilde
{H}=\left\{  0\right\}  $,\newline 3) $\widetilde{H}$ is an isomorphism from
$E$ to $\overline{E^{\ast}}$,\newline 4) $\operatorname*{Mat}\left(
H;e\right)  $ is non singular.
\end{theorem}

\textbf{Proof }: From $E^{\perp}=\left\{  x\in E\,/\,\forall y\in E\quad
H\left(  x,y\right)  =0\right\}  =\operatorname*{Ker}\widetilde{H}$ we
conclude that 1) is equivalent to~2). Since $\dim E=\dim\overline{E^{\ast}}$,
the condition $\operatorname*{Ker}\widetilde{H}=\left\{  0\right\}  $ means
that $\widetilde{H}$ is an isomorphism from $E$ to $\overline{E^{\ast}}$, and
it suffices to observe that $\operatorname*{Mat}\left(  H;e\right)  $ is the
matrix of $\widetilde{H}$ to complete the proof.%
\endproof

\begin{theorem}
Theorem~\ref{orth} can be improved if $H$ is non degenerate. We get:
\[
\dim F+\dim F^{\bot}=n\text{,\quad}F=\left(  F^{\bot}\right)  ^{\bot}\text{
and }\left(  F\cap G\right)  ^{\bot}=F^{\bot}+G^{\bot}.
\]
\end{theorem}

\textbf{Proof }: Let $\left(  e_{1},...,e_{p}\right)  $ denotes a basis of
$F$. As $\widetilde{H}$ is an isomorphism, the orthogonal
\[
F^{\bot}=\left\{  x\in E\,/\,\forall i\in\mathbb{N}_{p}\;H\left(
x,e_{i}\right)  =0\right\}  =\bigcap\limits_{i=1}^{p}\ker\widetilde{H}\left(
e_{i}\right)
\]
is the intersection of $p$ kernels of independant linear forms. Hence $\dim
F^{\bot}=n-p$. The inclusion $F\subset\left(  F^{\bot}\right)  ^{\bot}$ and
the equality $\dim\left(  F^{\bot}\right)  ^{\bot}=n-\dim F^{\bot}=\dim F$
give us the second result. It is sufficient to write the relation $\left(
F+G\right)  ^{\bot}=F^{\bot}\cap G^{\bot}$ of Theorem~\ref{orth} with
$F^{\bot}$ and $G^{\bot}$ instead of $F$ and $G$ to prove the last result.%
\endproof

\subsection{Isotropy}

\begin{definition}
Let $H$ be a hermitian form on$~E=\mathbb{F}_{t^{2}}^{N}$. A subspace $F$ of
$E$ is called \textbf{isotropic} if $F\cap F^{\bot}\neq\left\{  0\right\}  $.
A vector $x$ is \textbf{isotropic} if$~H\left(  x,x\right)  =0$.
\end{definition}

Note that a non null vector $x$ is isotropic if and only if the subspace
$\operatorname*{Vect}\left(  x\right)  $ generated by $x$ is isotropic.

\begin{theorem}
Let $H$ denote a hermitian form on$~E=\mathbb{F}_{t^{2}}^{N}$, and $F$ a
subspace of $E$. The following conditions are equivalent:\newline (i) The
restriction $H|_{F\times F}$ of $H$ to $F$ is non degenerate,\newline (ii) $F$
is non isotropic,\newline (iii) $E=F\oplus F^{\bot}$.
\end{theorem}

\textbf{Preuve }: Equivalence between $\left(  i\right)  $ and $\left(
ii\right)  $ follows from:
\begin{align*}
\left(  i\right)   &  \Leftrightarrow\forall x\in F\quad\left(  \left(
\forall y\in F\quad H\left(  x,y\right)  =0\right)  \Rightarrow x=0\right) \\
&  \Leftrightarrow\forall x\in F\quad\left(  x\in F^{\bot}\Rightarrow
x=0\right)  \Leftrightarrow F\cap F^{\bot}=\left\{  0\right\}  \Leftrightarrow
\left(  ii\right)  .
\end{align*}
We see at once that $\left(  iii\right)  $ implies $\left(  ii\right)  $. Let
us show that $\left(  i\right)  $ implies $\left(  iii\right)  $. If
$H|_{F\times F} $ is non degenerate, we already have $F\cap F^{\bot}=\left\{
0\right\}  $, and it only remains to prove that $E=F+F^{\bot}$. Let $x\in E$.
Let $l$ be the linear form $F\rightarrow\mathbb{F}_{t^{2}}\,;\,z\mapsto
H\left(  x,z\right)  $. Since $H|_{F\times F}$ is non degenerate, we can find
$y\in F$ such that $l=H|_{F\times F}\left(  y,.\right)  $, and it shows that
\[
\forall z\in F\quad l\left(  z\right)  =H\left(  x,z\right)  =H\left(
y,z\right)  .
\]
Thus $H\left(  x-y,z\right)  =0$ for all $z$ in $F$, and we conclude that
\[
\forall x\in E\quad\exists y\in F\quad x=\left(  x-y\right)  +y\text{ et
}x-y\in F^{\bot}.%
\endproof
\]

\section{Quadratic hermitian forms on $\mathbb{F}_{t^{2}}^{N}$\label{990924a}}

\begin{definition}
If $H$ denotes a hermitian form, the application
\[%
\begin{array}
[c]{cccc}%
q: & E & \longrightarrow & \mathbb{F}_{t}\\
& x & \longmapsto &  q\left(  x\right)  =H\left(  x,x\right)
\end{array}
\]
is called the \textbf{quadratic hermitian form }on $E$ associated to~$H$. We
denote by$~\operatorname*{QH}\left(  \mathbb{F}_{t^{2}}^{N}\right)  $ the
space of all quadratic hermitian forms on~$E$.
\end{definition}

With this Definition:

$\quad\left(  1\right)  \;\forall\lambda\in\mathbb{F}_{t^{2}}\quad\forall x\in
E\quad q\left(  \lambda x\right)  =\lambda^{t+1}q\left(  x\right)  =N\left(
\lambda\right)  q\left(  x\right)  ,$

$\quad\left(  1^{\prime}\right)  \;\forall\lambda\in\mathbb{F}_{t}\quad\forall
x\in E\quad q\left(  \lambda x\right)  =\lambda^{2}q\left(  x\right)  ,$

$\quad\left(  2\right)  \;\forall x,y\in E\;H\left(  x,y\right)  =\frac{1}%
{4}\left[  q\left(  x+y\right)  -q\left(  x-y\right)  \right]  +\frac
{1}{2\left(  \alpha-\overline{\alpha}\right)  }\left[  q\left(  x+\alpha
y\right)  -q\left(  x+\overline{\alpha}y\right)  \right]  $.

Result $\left(  1^{\prime}\right)  $ explain the name ''quadratic'' when we
restrict our attention on$~\mathbb{F}_{t}$. Result $\left(  2\right)  $ is
true when $t$ is~odd (see Lemma~\ref{990921a}). From now on we assume that~$t
$ is~odd.

\begin{theorem}
Let $\mathbb{F}_{t}^{E}$ denotes the $\mathbb{F}_{t}$-vector space of all
applications from $E$ to $\mathbb{F}_{t}$. The function
\[%
\begin{array}
[c]{cccc}%
\Psi: & \operatorname*{H}\left(  \mathbb{F}_{t^{2}}^{N}\right)  &
\longrightarrow & \mathbb{F}_{t}^{E}\\
& H & \longmapsto &  q\text{ such that }q\left(  x\right)  =H\left(
x,x\right)
\end{array}
\]
is $\mathbb{F}_{t}$-linear and one to one. We have $\operatorname*{Im}%
\Psi=\operatorname*{QH}\left(  \mathbb{F}_{t^{2}}^{N}\right)  $ and $\Psi$
induces an isomorphism from $\operatorname*{H}\left(  \mathbb{F}_{t^{2}}%
^{N}\right)  $ onto $\operatorname*{QH}\left(  \mathbb{F}_{t^{2}}^{N}\right)
$. Hence$~\dim_{\mathbb{F}_{t}}\operatorname*{QH}\left(  \mathbb{F}_{t^{2}%
}^{N}\right)  =N^{2}$.
\end{theorem}

\textbf{Proof }: Result $\left(  2\right)  $ shows that$~\Psi$ is one to one.%
\endproof

\begin{definition}
Let $q\in\operatorname*{QH}\left(  \mathbb{F}_{t^{2}}^{N}\right)  $. The
unique hermitian form $H$ satisfying $\Psi\left(  H\right)  =q$ is called the
\textbf{polar form} of$~q$, and is given by result$~\left(  2\right)  $. The
isomorphism $\operatorname*{QH}\left(  \mathbb{F}_{t^{2}}^{N}\right)
\simeq\operatorname*{H}\left(  \mathbb{F}_{t^{2}}^{N}\right)  $ allows us to
construct the same objects from a quadratic hermitian form or from a hermitian
form. For instance, the \textbf{kernel} and the \textbf{rank} of$~q$ will be
those of the associated polar form.
\end{definition}

Let $M=\left(  a_{ij}\right)  $ denotes the matrix of a hermitian form $H$ in
a basis of~$E$. Then
\begin{align*}
H\left(  x,x\right)   &  =\sum\limits_{i,j}a_{ij}\overline{x}_{i}x_{j}%
=\sum\limits_{i=1}^{N}a_{ii}x_{i}^{t+1}+\sum\limits_{1\leq i<j\leq N}\left(
a_{ij}\overline{x}_{i}x_{j}+a_{ji}\overline{x}_{j}x_{i}\right) \\
&  =\sum\limits_{i=1}^{N}a_{ii}x_{i}^{t+1}+\sum\limits_{1\leq i<j\leq
N}\left(  a_{ij}\overline{x}_{i}x_{j}+\left(  a_{ij}\overline{x}_{i}%
x_{j}\right)  ^{t}\right)  .
\end{align*}
We can say that a quadratic hermitian form $q$\ on $E$ is an application from
$E$ to $\mathbb{F}_{t}$ defined by
\[
\forall x\in E\quad q\left(  x\right)  =\sum\limits_{i=1}^{N}a_{ii}%
\operatorname*{N}\nolimits_{\mathbb{F}_{t^{2}}/\mathbb{F}_{t}}\left(
x_{i}\right)  +\sum\limits_{1\leq i<j\leq N}\operatorname*{Tr}%
\nolimits_{\mathbb{F}_{t^{2}}/\mathbb{F}_{t}}\left(  a_{ij}\overline{x}%
_{i}x_{j}\right)
\]
where $\left(  x_{1},...,x_{N}\right)  $ are the coordinates of $x$ in a
basis, $a_{ii}\in\mathbb{F}_{t}$ and $a_{ij}\in\mathbb{F}_{t^{2}}$ for all
$i\neq j$. In fact Theorem~\ref{990922a} ensures us the existence of a
$H$-orthogonal basis of $E$. In such a basis $q\left(  x\right)  =\sum
_{i=1}^{r}x_{i}^{t+1}$ for all$~x\in E$.

Next Theorem provides another criterion for~$q$:

\begin{theorem}
If $q\in\mathbb{F}_{t}^{E}$ satisfies $\left(  1\right)  $ and if $H$ defined
by $\left(  2\right)  $ is sesquilinear, then$~q$ is the quadratric hermitian
form associated with the hermitian form$~H$.
\end{theorem}

\textbf{Proof }: $\left(  1\right)  $ and $\left(  2\right)  $ show that
\begin{align*}
H\left(  x,x\right)   &  =\frac{1}{4}\left[  q\left(  2x\right)  -q\left(
0\right)  \right]  +\frac{1}{2\left(  \alpha-\overline{\alpha}\right)
}\left[  q\left(  \left(  1+\alpha\right)  x\right)  -q\left(  \left(
1+\overline{\alpha}\right)  x\right)  \right] \\
&  =\frac{1}{4}\left[  2^{t+1}q\left(  x\right)  \right]  +\frac{1}{2\left(
\alpha-\overline{\alpha}\right)  }\left[  \left(  1+\alpha\right)
^{t+1}-\left(  1+\overline{\alpha}\right)  ^{t+1}\right]  q\left(  x\right) \\
&  =q\left(  x\right)  +\frac{1}{2\left(  \alpha-\overline{\alpha}\right)
}\left[  \left(  1+\alpha^{t}\right)  \left(  1+\alpha\right)  -\left(
1+\overline{\alpha}^{t}\right)  \left(  1+\overline{\alpha}\right)  \right]
q\left(  x\right)  =q\left(  x\right)  .
\end{align*}
By hypothesis, the form $H$ is sesquilinear, and it remains to prove the
hermitian symetry. We have
\[
H\left(  y,x\right)  =\frac{1}{4}\left[  q\left(  y+x\right)  -q\left(
y-x\right)  \right]  +\frac{1}{2\left(  \alpha-\overline{\alpha}\right)
}\left[  q\left(  y+\alpha x\right)  -q\left(  y+\overline{\alpha}x\right)
\right]  .
\]
Condition 2) of Lemma~\ref{990921a} yields
\[
q\left(  y+\alpha x\right)  -q\left(  y+\overline{\alpha}x\right)  =q\left(
x+\overline{\alpha}y\right)  -q\left(  x+\alpha y\right)
\]
hence
\[
H\left(  y,x\right)  =\frac{1}{4}\left[  q\left(  x+y\right)  -q\left(
x-y\right)  \right]  -\frac{1}{2\left(  \alpha-\overline{\alpha}\right)
}\left[  q\left(  x+\alpha y\right)  -q\left(  x+\overline{\alpha}y\right)
\right]  .
\]
As $q$\ takes its values in $\mathbb{F}_{t}$, we conclude that$~H\left(
y,x\right)  =\overline{H\left(  x,y\right)  }$.~%
\endproof
$\medskip$

\textbf{Remark }: The proof above also gives that a sesquilinear form $H$ is
hermitian if and only if $H\left(  x,x\right)  \in\mathbb{F}_{t}$ for
all$~x\in E$.

\section{Equivalence between quadratic hermitian forms}

\begin{definition}
Two hermitian forms (resp. quadratic hermitian forms) $\varphi_{1}$ and
$\varphi_{2}$ (resp. $q_{1}$ and $q_{2}$) are called \textbf{equivalent}, and
we note $\varphi_{1}\sim\varphi_{2}$ (resp. $q_{1}\sim q_{2}$), if there
exists an automorphism $u$ of $E$ such that
\[
\forall x,y\in E\quad\varphi_{2}\left(  x,y\right)  =\varphi_{1}\left(
u\left(  x\right)  ,u\left(  y\right)  \right)  \quad\text{(resp.\ }\forall
x\in E\quad q_{2}\left(  x\right)  =q_{1}\left(  u\left(  x\right)  \right)
\text{).}%
\]
\end{definition}

\begin{theorem}
Let $q_{1}$ and $q_{2}$ denote two quadratic hermitian forms whose polar forms
are $\varphi_{1}$ and $\varphi_{2}$. The following conditions are
equivalent:\newline i) $q_{1}$ and $q_{2}$ are equivalent,\newline ii) $q_{1}$
and $q_{2}$ have the same matrix but in different basis,\newline iii) $q_{1}$
and $q_{2}$ have same rank.
\end{theorem}

\textbf{Proof }: i) $\Leftrightarrow$ ii) : Let $e=\left(  e_{1}%
,...,e_{N}\right)  $ denote a basis of $E$. We have
\begin{align*}
\left(  \varphi_{1}\sim\varphi_{2}\right)   &  \Leftrightarrow\exists
u\in\operatorname*{GL}\left(  E\right)  \quad\forall x,y\in E\quad\varphi
_{2}\left(  x,y\right)  =\varphi_{1}\left(  u\left(  x\right)  ,u\left(
y\right)  \right) \\
&  \Leftrightarrow\exists u\in\operatorname*{GL}\left(  E\right)  \quad\forall
i,j\in\mathbb{N}_{n}\quad\varphi_{2}\left(  e_{i},e_{j}\right)  =\varphi
_{1}\left(  u\left(  e_{i}\right)  ,u\left(  e_{j}\right)  \right) \\
&  \Leftrightarrow\exists u\in\operatorname*{GL}\left(  E\right)
\quad\operatorname*{Mat}\left(  \varphi_{2};e\right)  =\operatorname*{Mat}%
\left(  \varphi_{1};u\left(  e\right)  \right)  ,
\end{align*}
therefore i) implies ii).

Conversely, if there are two basis $e$ and $e^{\prime}$ such that
$\operatorname*{Mat}\left(  \varphi_{2};e\right)  =\operatorname*{Mat}\left(
\varphi_{1};e^{\prime}\right)  $, we can define an automorphism $u$ in $E$
with $u\left(  e\right)  =e^{\prime}$, and use the above equivalences to
obtain$~\varphi_{1}\sim\varphi_{2}$.\smallskip

\quad ii) $\Leftrightarrow$ iii) : If $\operatorname*{Mat}\left(  \varphi
_{2};e\right)  =\operatorname*{Mat}\left(  \varphi_{1};e^{\prime}\right)  $
then $\varphi_{1}$ and $\varphi_{2}$ have same rank. Conversely, if
$\varphi_{1} $ and $\varphi_{2}$ have same rank $r$, Theorem~\ref{990922a}
provides two basis $e$ and $e^{\prime}$ such that $\operatorname*{Mat}\left(
\varphi_{2};e\right)  $ and $\operatorname*{Mat}\left(  \varphi_{1};e^{\prime
}\right)  $ are\ equal to the diagonal matrix $\operatorname*{Diag}\left(
1,...,1,0,...,0\right)  $ with $r$ numbers~$1$.~%
\endproof

\section{Quadratic hermitian forms on $\mathbb{F}_{t}^{2N}$}

Suppose $t$ odd. Let $H:\mathbb{F}_{t^{2}}^{N}\times\mathbb{F}_{t^{2}}%
^{N}\rightarrow\mathbb{F}_{t^{2}}$ be a sesquilinear form and $\alpha$ denotes
an element of $\mathbb{F}_{t^{2}}$ with $\mathbb{F}_{t^{2}}=\mathbb{F}%
_{t}\left(  \alpha\right)  $. The application $\iota:\mathbb{F}_{t}%
^{2N}\rightarrow\mathbb{F}_{t^{2}}^{N}$ defined by
\[
\iota\left(  x_{1},...,x_{2N}\right)  =\left(  x_{1}+\alpha x_{2}%
,...,x_{2N-1}+\alpha x_{2N}\right)
\]
is an $\mathbb{F}_{t}$-vector space isomorphism. Since $H:\mathbb{F}_{t^{2}%
}^{N}\times\mathbb{F}_{t^{2}}^{N}\rightarrow\mathbb{F}_{t^{2}}$ is
$\mathbb{F}_{t}$-bilinear, it will be the same with $\mathbb{F}_{t}^{2N}%
\times\mathbb{F}_{t}^{2N}\rightarrow\mathbb{F}_{t^{2}};\left(  x,y\right)
\mapsto H\left(  \iota x,\iota y\right)  $. Roughly speaking, we want to work
with functions with values in$~\mathbb{F}_{t}$, thus it is convenient to define:

\begin{definition}
The \textbf{quadratic hermitian form }$f$\textbf{\ on} $\mathbb{F}_{t}^{2N}%
$\textbf{\ associated with }$H$ is
\[%
\begin{array}
[c]{cccc}%
f: & \mathbb{F}_{t}^{2N} & \rightarrow & \mathbb{F}_{t}\\
& x & \mapsto &  H\left(  \iota x,\iota x\right)  .
\end{array}
\]
We denote by $\operatorname*{QH}\left(  \mathbb{F}_{t}^{2N}\right)  $ the
vector space of all quadratic hermitian forms $f$\ on $\mathbb{F}_{t}^{2N}$.
\end{definition}

It is clear that the function
\[%
\begin{array}
[c]{ccc}%
\operatorname*{QH}\left(  \mathbb{F}{_{t^{2}}^{N}}\right)  & \rightarrow &
\operatorname*{QH}\left(  \mathbb{F}_{t}^{2N}\right) \\
q & \mapsto &  f
\end{array}
\]
where $f\left(  x\right)  =H\left(  \iota x,\iota x\right)  $ when $q\left(
x\right)  =H\left(  x,x\right)  $, is an isomorphism between $\mathbb{F}_{t}%
$-vector spaces, hence $\dim_{\mathbb{F}_{t}}\operatorname*{QH}\left(
\mathbb{F}_{t}^{2N}\right)  =N^{2}$.

\begin{theorem}
\label{991008a}Suppose $t$ odd. The quadratic hermitian form $f$\ on
$\mathbb{F}_{t}^{2N}$\ associated with $H$ is a $\mathbb{F}_{t}$-quadratic
form associated with the bilinear form $\frac{1}{2}B$, where
\[%
\begin{array}
[c]{cccc}%
B: & \mathbb{F}_{t}^{2N}\times\mathbb{F}_{t}^{2N} & \rightarrow &
\mathbb{F}_{t}\\
& \left(  x,y\right)  & \mapsto &  f\left(  x+y\right)  -f\left(  x\right)
-f\left(  y\right)  .
\end{array}
\]
We have $B\left(  x,y\right)  =H\left(  \iota x,\iota y\right)  +H\left(
\iota x,\iota y\right)  ^{t}=\operatorname*{Tr}\nolimits_{\mathbb{F}_{t^{2}%
}/\mathbb{F}_{t}}\left(  H\left(  \iota x,\iota y\right)  \right)  .$
\end{theorem}

\textbf{Proof }: It is a simple matter to see that $f$ is a $\mathbb{F}_{t}%
$-quadratic form because $f\left(  x\right)  $ is a homogeneous polynomial of
degree $2$ in the coodinates of $x$ and with coefficients in$~\mathbb{F}_{t}$.
Indeed, it suffices to use a $H$-orthogonal basis of $\mathbb{F}_{t^{2}}^{N}$
to get $q\left(  x\right)  =\sum_{i=1}^{r}x_{i}^{t+1}$ for all $x\in E$
(Theorem~\ref{990922a}) and
\[
f\left(  x\right)  =q\left(  \iota x\right)  =\sum\limits_{i=1}^{r}\left(
u_{i}+\alpha v_{i}\right)  ^{t+1}=\sum\limits_{i=1}^{r}u_{i}^{2}+\alpha
^{t+1}v_{i}^{2}+\left(  \alpha+\alpha^{t}\right)  u_{i}v_{i}%
\]
with $x=\left(  u_{1},v_{1},...,u_{N},v_{N}\right)  \in\mathbb{F}_{t}^{2N}$
and $\alpha^{t+1}\in\mathbb{F}_{t}$.

Then it is easy to check that $f\left(  x+y\right)  =f\left(  x\right)
+f\left(  y\right)  +H\left(  \iota x,\iota y\right)  +H\left(  \iota x,\iota
y\right)  ^{t} $.%
\endproof

\begin{definition}
For simplicity, we also say that $B$ is a bilinear form associated with~$f$.
\end{definition}

\textbf{Remark :} If $t$ is even, say $t=2^{\zeta}$, a symmetric bilinear form
on$~\mathbb{F}_{t}^{2N}$ is
\[
B\left(  x,y\right)  =\sum\limits_{i}a_{ii}x_{i}y_{i}+\sum\limits_{i<j}%
a_{ij}\left(  x_{i}y_{j}+x_{j}y_{i}\right)  \quad\quad\left(  \ast\right)
\]
and the quadratic form associated to $B$ is $f\left(  x\right)  =B\left(
x,x\right)  =\sum_{i}a_{ii}x_{i}^{2}$. A quadratic form on$~\mathbb{F}%
_{2^{\zeta}}^{2N}$ will be an homogeneous polynomial of degree~$2$ in the
coordinates $x_{1},...,x_{2N}$ with no diagonal term$~a_{ij}x_{i}x_{j}$.
Conversely, if $f\left(  x\right)  =\sum_{i}a_{ii}x_{i}^{2}$ is a quadratic
form on$~\mathbb{F}_{t}^{2N}$, there exists an infinity of symmetric bilinear
forms $B\left(  x,y\right)  $ such that $B\left(  x,x\right)  =f\left(
x\right)  $, and this is different from the usual case. Indeed, it suffices to
choose any coefficients $a_{ij}$ ( $i<j$) in $\mathbb{F}_{2^{\zeta}}$ and to
define $B\left(  x,y\right)  $ by $\left(  \ast\right)  $ to get $B\left(
x,x\right)  =f\left(  x\right)  $. In this case, the map $B\left(  x,y\right)
=f\left(  x+y\right)  -f\left(  x\right)  -f\left(  y\right)  $ of
Theorem~\ref{991008a} will never be the bilinear form associated with~$f$ as
$B\left(  x,x\right)  =f\left(  2x\right)  -2f\left(  x\right)  =0$.\medskip

The kernel of $B$ is
\[
\operatorname*{Ker}B=\left\{  x\in\mathbb{F}_{t}^{2N}\,/\,\forall
y\in\mathbb{F}_{t}^{2N}\;B\left(  x,y\right)  =0\right\}  ,
\]
and the orthogonal of $\operatorname*{Ker}B$ for the usual inner product in
$\mathbb{F}_{t}^{2N}$ is
\[
\left(  \operatorname*{Ker}B\right)  ^{\bot}=\left\{  x\in\mathbb{F}_{t}%
^{2N}\,/\,\forall y\in\operatorname*{Ker}B\;x.y=x_{1}.y_{1}+...+x_{2N}%
.y_{2N}=0\right\}  .
\]
Since the usual inner product $x.y$ is only a non degenerate bilinear form on
$\mathbb{F}_{t}^{2N}$, we have $\dim\operatorname*{Ker}B+\dim\left(
\operatorname*{Ker}B\right)  ^{\bot}=2N$ but we can't say that $\mathbb{F}%
_{t}^{2N}=\operatorname*{Ker}B\oplus\left(  \operatorname*{Ker}B\right)
^{\bot}$. With these notations~:

\begin{theorem}
\label{990925b}We have $\iota\left(  \operatorname*{Ker}B\right)
=\operatorname*{Ker}H$. Thus $\iota$ induces a $\mathbb{F}_{t}$-isomorphism
from $\operatorname*{Ker}B$ onto $\operatorname*{Ker}H$ and
$\operatorname*{rk}f=\operatorname*{rk}B=2\operatorname*{rk}H$.
\end{theorem}

\textbf{Proof }: Let $\psi$ denotes a non trivial additive character
on~$\mathbb{F}_{t}$. The map $\psi^{\prime}=\psi\circ\operatorname*{Tr}%
\nolimits_{\mathbb{F}_{t^{2}}/\mathbb{F}_{t}}$ is a non trivial additive
character on~$\mathbb{F}_{t^{2}}$ and Lemma~\ref{990925a} gives:
\begin{align*}
\left(  x\in\operatorname*{Ker}B\right)   &  \Leftrightarrow\forall
y\in\mathbb{F}_{t}^{2N}\quad B\left(  x,y\right)  =\operatorname*{Tr}%
\nolimits_{\mathbb{F}_{t^{2}}/\mathbb{F}_{t}}\left(  H\left(  \iota x,\iota
y\right)  \right)  =0\\
&  \Leftrightarrow\sum\limits_{y\in\mathbb{F}_{t}^{2N}}\psi\left(
\operatorname*{Tr}\nolimits_{\mathbb{F}_{t^{2}}/\mathbb{F}_{t}}\left(
H\left(  \iota x,\iota y\right)  \right)  \right)  \neq0\Leftrightarrow
\sum\limits_{z\in\mathbb{F}_{t^{2}}^{N}}\psi\left(  \operatorname*{Tr}%
\nolimits_{\mathbb{F}_{t^{2}}/\mathbb{F}_{t}}\left(  H\left(  \iota
x,z\right)  \right)  \right)  \neq0\\
&  \Leftrightarrow\sum\limits_{z\in\mathbb{F}_{t^{2}}^{N}}\psi^{\prime}\left(
H\left(  \iota x,z\right)  \right)  \neq0\Leftrightarrow\forall z\in
\mathbb{F}_{t^{2}}^{N}\quad H\left(  \iota x,z\right)  =0\Leftrightarrow\iota
x\in\operatorname*{Ker}H\text{.}%
\end{align*}
Hence $\iota\left(  \operatorname*{Ker}B\right)  \subset\operatorname*{Ker}H$.
Since $\iota$ is a $\mathbb{F}_{t}$-isomorphism, the above equivalences imply
the inverse inclusion. To complete the proof, we write
\[
\operatorname*{rk}f=\operatorname*{rk}B=2N-\dim_{\mathbb{F}_{t}}%
\operatorname*{Ker}B=2N-2\dim_{\mathbb{F}_{t^{2}}}\operatorname*{Ker}%
H=2\operatorname*{rk}H.\text{%
\endproof
}%
\]

\begin{theorem}
\label{991008d}1) There is an endomorphism $T$ of $\mathbb{F}_{t}^{2N}$ such
that $B\left(  x,y\right)  =T\left(  x\right)  .y$ for all $(x,y)\in
\mathbb{F}_{t}^{2N}\times\mathbb{F}_{t}^{2N}$.\newline 2) We have
$\operatorname*{Ker}T=\operatorname*{Ker}B$,\ $\operatorname*{Im}T=\left(
\operatorname*{Ker}B\right)  ^{\bot}$\ and $\operatorname*{Ker}T\subset
f^{-1}\left(  0\right)  $.
\end{theorem}

\textbf{Proof }: 1) Since the inner product is non degenerate, for all
$x\in\mathbb{F}_{t}^{2N}$ we can find $T\left(  x\right)  \in\mathbb{F}%
_{t}^{2N}$ such that $B\left(  x,y\right)  =T\left(  x\right)  .y$ for all
$y\in\mathbb{F}_{t}^{2N}$. From $B\left(  \lambda x+x^{\prime},y\right)
=\lambda B\left(  x,y\right)  +B\left(  x^{\prime},y\right)  $ we deduce
$\left[  T\left(  \lambda x+x^{\prime}\right)  -\lambda T\left(  x\right)
-T\left(  x^{\prime}\right)  \right]  .y=0$ for all $y\in\mathbb{F}_{t}^{2N}$,
hence
\[
T\left(  \lambda x+x^{\prime}\right)  -\lambda T\left(  x\right)  -T\left(
x^{\prime}\right)  =0
\]
and the linearity of$~T$ follows.\smallskip

2) The first equality is a consequence of
\[
x\in\operatorname*{Ker}T\Leftrightarrow\left(  \forall y\in\mathbb{F}_{t}%
^{2N}\quad T\left(  x\right)  .y=0\right)  \Leftrightarrow\left(  \forall
y\in\mathbb{F}_{t}^{2N}\quad B\left(  x,y\right)  =0\right)  \Leftrightarrow
x\in\operatorname*{Ker}B.
\]
If $z\in\mathbb{F}_{t}^{2N}$ and if $u\in\operatorname*{Ker}B$, then $T\left(
z\right)  .u=B\left(  z,u\right)  =0$, hence $\operatorname*{Im}%
T\subset\left(  \operatorname*{Ker}B\right)  ^{\bot}$. This inclusion is an
equality because
\[
\dim\left(  \operatorname*{Im}T\right)  =2N-\dim\left(  \operatorname*{Ker}%
T\right)  =2N-\dim\left(  \operatorname*{Ker}B\right)  =\dim\left(  \left(
\operatorname*{Ker}B\right)  ^{\bot}\right)  .
\]
If $x\in\operatorname*{Ker}T=\operatorname*{Ker}B$ then $f\left(  x\right)
=H\left(  \iota x,\iota x\right)  =0$ from Theorem~\ref{990925b}, thus
$\operatorname*{Ker}T\subset f^{-1}\left(  0\right)  $.%
\endproof

\section{Exponential sums S$\left(  f,v\right)  $}

Let us denote by $\psi$ the additive character on $\mathbb{F}_{t}$ defined by
\[
\psi\left(  x\right)  =\exp\left(  \frac{i2\pi}{p}\operatorname*{Tr}%
\nolimits_{\mathbb{F}_{t}/\mathbb{F}_{p}}\left(  x\right)  \right)  \text{.}%
\]
If $v\in\mathbb{F}_{t}^{2N}$, we consider the exponential sum associated to
$f$ and $v$:
\[
S\left(  f,v\right)  =\sum\limits_{x\in\mathbb{F}_{t}^{2N}}\psi\left(
f\left(  x\right)  +v.x\right)  .
\]

\begin{lemma}
\label{990925a} Let $\psi$ denotes a non trivial additive character
on~$\mathbb{F}_{t}$, $V$ a $\mathbb{F}_{t}$-vector space of finite dimension
$m$, and $l:V\rightarrow\mathbb{F}_{t}$ a linear form on $V$. Then
\[
\sum\limits_{y\in V}\psi\left(  l\left(  y\right)  \right)  =\left\{
\begin{array}
[c]{ll}%
t^{m} & \;\text{if }l=0,\\
0 & \;\text{if }l\neq0.
\end{array}
\right.
\]
\end{lemma}

\textbf{Proof }: The map $\psi\circ l$ is an additive character on
$V\simeq\mathbb{F}_{t}^{m}$ and we can apply the orthogonality relation
(Theorem~\ref{991008b}).%
\endproof

\begin{lemma}
\label{991008c}
\[
\sum\limits_{x\in\mathbb{F}_{t^{m}}}\psi\left(  \operatorname*{N}%
\nolimits_{\mathbb{F}_{t^{m}}/\mathbb{F}_{t}}\left(  x\right)  \right)
=\frac{t-t^{m}}{t-1}\text{.}%
\]
\end{lemma}

\textbf{Proof }: From Lemma \ref{990921b} it follows that $N_{\mathbb{F}%
_{t^{m}}/\mathbb{F}_{t}}:\mathbb{F}_{t^{m}}^{\ast}\rightarrow\mathbb{F}%
_{t}^{\ast}$ \ is a multiplicative group epimorphism and that $\left|
N_{\mathbb{F}_{t^{m}}/\mathbb{F}_{t}}{}^{-1}\left(  b\right)  \right|
=\frac{t^{m}-1}{t-1}$ for all$~b\in\mathbb{F}_{t}^{\ast}$. Hence
\[
\sum\limits_{x\in\mathbb{F}_{t^{m}}}\psi\left(  \operatorname*{N}%
\nolimits_{\mathbb{F}_{t^{m}}/\mathbb{F}_{t}}\left(  x\right)  \right)
=1+\sum\limits_{x\in\mathbb{F}_{t^{m}}^{\ast}}\psi\left(  \operatorname*{N}%
\nolimits_{\mathbb{F}_{t^{m}}/\mathbb{F}_{t}}\left(  x\right)  \right)
=1+\frac{t^{m}-1}{t-1}\sum\limits_{z\in\mathbb{F}_{t}^{\ast}}\psi\left(
z\right)  .
\]
The use of the orthogonality relation $\sum_{z\in\mathbb{F}_{t}^{\ast}}%
\psi\left(  z\right)  =-1$ completes the proof.%
\endproof
\medskip

We are now ready to give another proof of the main result in~\cite{cherdieu}.
In fact, a small mistake occured in Proposition~3 of \cite{cherdieu} as
$A\left(  s,v\right)  $ do not depends on $f\left(  u\right)  $ but on
$\operatorname*{Tr}\nolimits_{\mathbb{F}_{t}/\mathbb{F}_{s}}\left(  f\left(
u\right)  \right)  $, as we shall see below.

\begin{theorem}
\label{cher}\{\cite{cherdieu}, Th.~2 and Prop.~3\} Let $v\in\mathbb{F}%
_{t}^{2N}$ and let $f$ denote a quadratic hermitian form of rank~$2\rho$
in$~\mathbb{F}_{t}^{2N}$. Consider the extensions $\mathbb{F}_{p}%
\subset\mathbb{F}_{s}\subset\mathbb{F}_{t}\subset\mathbb{F}_{t^{2}}$ and let
$a\in\mathbb{F}_{s}^{\ast}$.\newline 1) If $v\in\left(  \operatorname*{Ker}%
B\right)  ^{\bot}=\operatorname*{Im}T$, we can find $u\in\mathbb{F}_{t}^{2N}$
such that $v=T\left(  u\right)  $. Then
\[
S\left(  af,v\right)  =\left(  -1\right)  ^{\rho}t^{2N-\rho}\psi\left(
-a^{-1}f\left(  u\right)  \right)
\]
and $\sum_{a\in\mathbb{F}_{s}^{\ast}}S\left(  af,v\right)  =\left(  -1\right)
^{\rho}t^{2N-\rho}A\left(  s,v\right)  $ where
\[
A\left(  s,v\right)  =\left\{
\begin{array}
[c]{ll}%
s-1 & \text{if }\operatorname*{Tr}\nolimits_{\mathbb{F}_{t}/\mathbb{F}_{s}%
}\left(  f\left(  u\right)  \right)  =0,\\
-1 & \text{else.}%
\end{array}
\right.
\]
2) If $v\notin\left(  \operatorname*{Ker}B\right)  ^{\bot}$ then $S\left(
af,v\right)  =0$.
\end{theorem}

\textbf{Proof }: Without loss of generality, we can assume that $f$ is given
in the standard form $f\left(  x\right)  =H\left(  y,y\right)  =q\left(
y\right)  =y_{1}^{t+1}+...+y_{\rho}^{t+1}$ where $y=\iota\left(  x\right)
\in\mathbb{F}_{t^{2}}^{N}$.\smallskip

\quad1) $\alpha)$ We first compute $S\left(  f,v\right)  $. Since $v=T\left(
u\right)  $,
\[
f\left(  x\right)  +v.x=f\left(  x\right)  +T\left(  u\right)  .x=f\left(
x\right)  +B\left(  u,x\right)  =f\left(  u+x\right)  -f\left(  u\right)
\]
and
\[
S\left(  f,v\right)  =\sum\limits_{x\in\mathbb{F}_{t}^{2N}}\psi\left(
f\left(  x\right)  +v.x\right)  =\sum\limits_{x\in\mathbb{F}_{t}^{2N}}%
\psi\left(  f\left(  u+x\right)  -f\left(  u\right)  \right)  .
\]
Define $z=\iota u$. Then
\[
f\left(  u+x\right)  -f\left(  u\right)  =q\left(  z+y\right)  -q\left(
z\right)  =\sum\limits_{k=1}^{\rho}\left[  \left(  z_{k}+y_{k}\right)
^{t+1}-z_{k}^{t+1}\right]
\]
and
\begin{align*}
S\left(  f,v\right)   &  =\sum\limits_{y_{1},...,y_{N}\in\mathbb{F}_{t^{2}}%
}\prod\limits_{k=1}^{\rho}\psi\left(  \left(  z_{k}+y_{k}\right)  ^{t+1}%
-z_{k}^{t+1}\right) \\
&  =t^{2\left(  N-\rho\right)  }\sum\limits_{y_{1},...,y_{\rho}\in
\mathbb{F}_{t^{2}}}\prod\limits_{k=1}^{\rho}\psi\left(  \left(  z_{k}%
+y_{k}\right)  ^{t+1}-z_{k}^{t+1}\right) \\
&  =t^{2\left(  N-\rho\right)  }\xi\prod\limits_{k=1}^{\rho}\psi\left(
-z_{k}^{t+1}\right)
\end{align*}
where $\xi=\sum\limits_{y_{1},...,y_{\rho}\in\mathbb{F}_{t^{2}}}%
\prod\limits_{k=1}^{\rho}\psi\left(  \left(  z_{k}+y_{k}\right)
^{t+1}\right)  $. We have
\[
\xi=\sum\limits_{y_{1},...,y_{\rho-1}\in\mathbb{F}_{t^{2}}}\left(
\prod\limits_{k=1}^{\rho-1}\psi\left(  \left(  z_{k}+y_{k}\right)
^{t+1}\right)  \right)  \left(  \sum\limits_{y_{\rho}\in\mathbb{F}_{t^{2}}%
}\psi\left(  \left(  z_{\rho}+y_{\rho}\right)  ^{t+1}\right)  \right)  .
\]
Lemma~\ref{991008c} gives $\sum\limits_{y_{\rho}\in\mathbb{F}_{t^{2}}}%
\psi\left(  \left(  z_{\rho}+y_{\rho}\right)  ^{t+1}\right)  =\sum
\limits_{y\in\mathbb{F}_{t^{2}}}\psi\left(  y^{t+1}\right)  =-t$, hence
\[
\xi=\left(  -t\right)  \sum\limits_{y_{1},...,y_{\rho-1}\in\mathbb{F}_{t^{2}}%
}\left(  \prod\limits_{k=1}^{\rho-1}\psi\left(  \left(  z_{k}+y_{k}\right)
^{t+1}\right)  \right)  .
\]
We proceed to obtain $\xi=\left(  -t\right)  ^{\rho}$, and so
\begin{align*}
S\left(  f,v\right)   &  =\left(  -1\right)  ^{\rho}t^{2N-\rho}\prod
\limits_{k=1}^{\rho}\psi\left(  -z_{k}^{t+1}\right)  =\left(  -1\right)
^{\rho}t^{2N-\rho}\psi\left(  -z_{1}^{t+1}-...-z_{\rho}^{t+1}\right) \\
&  =\left(  -1\right)  ^{\rho}t^{2N-\rho}\psi\left(  -q\left(  z\right)
\right)  .
\end{align*}
Since $q\left(  z\right)  =H\left(  \iota u,\iota u\right)  =f\left(
u\right)  $, we see that $S\left(  f,v\right)  =\left(  -1\right)  ^{\rho
}t^{2N-\rho}\psi\left(  -f\left(  u\right)  \right)  $.\medskip

$\quad\beta)$ Let us compute $S\left(  af,v\right)  $. By the above applied
with $f_{a}=af$ instead of $f$, we obtain $S\left(  af,v\right)  =\left(
-1\right)  ^{\rho}t^{2N-\rho}\psi\left(  -af\left(  u_{a}\right)  \right)  $
where $u_{a}$ satisfies $v=T_{a}u_{a}$ and $T_{a}$ is defined by
\[
T_{a}\left(  x\right)  .y=f_{a}\left(  x+y\right)  -f_{a}\left(  x\right)
-f_{a}\left(  y\right)  =a\left(  f\left(  x+y\right)  -f\left(  x\right)
-f\left(  y\right)  \right)  =a\left(  T\left(  x\right)  .y\right)  .
\]
Hence $T_{a}=aT$. We have $v=T_{a}u_{a}=aT\left(  u_{a}\right)  =T\left(
au_{a}\right)  $, and we can take $u=au_{a}$. This gives $S\left(
af,v\right)  =\left(  -1\right)  ^{\rho}t^{2N-\rho}\psi\left(  -af\left(
a^{-1}u\right)  \right)  =\left(  -1\right)  ^{\rho}t^{2N-\rho}\psi\left(
-a^{-1}f\left(  u\right)  \right)  $.\medskip

$\quad\gamma)$ By the above
\begin{align*}
\sum\limits_{a\in\mathbb{F}_{s}^{\ast}}S\left(  af,v\right)   &  =\left(
-1\right)  ^{\rho}t^{2N-\rho}\sum\limits_{a\in\mathbb{F}_{s}^{\ast}}%
\psi\left(  -a^{-1}f\left(  u\right)  \right) \\
&  =\left(  -1\right)  ^{\rho}t^{2N-\rho}\sum\limits_{a\in\mathbb{F}_{s}%
^{\ast}}\psi^{\prime}\left(  -a^{-1}\operatorname*{Tr}\nolimits_{\mathbb{F}%
_{t}/\mathbb{F}_{s}}\left(  f\left(  u\right)  \right)  \right)
\end{align*}
where $\psi^{\prime}$ is the additive character $\psi^{\prime}\left(
x\right)  =\exp\left(  \frac{i2\pi}{p}\operatorname*{Tr}\nolimits_{\mathbb{F}%
_{s}/\mathbb{F}_{p}}\left(  x\right)  \right)  $ on$~\mathbb{F}_{s}$. The map
$z\mapsto\psi^{\prime}\left(  cz\right)  $ describes the set of additive
characters on $\mathbb{F}_{s}$ when $c$ describes$~\mathbb{F}_{s}$,
consequently the orthogonality relation (Theorem~\ref{991008b}) yields
\begin{align*}
\sum\limits_{a\in\mathbb{F}_{s}^{\ast}}S\left(  af,v\right)   &  =\left(
-1\right)  ^{\rho}t^{2N-\rho}\left(  -1+\sum\limits_{\chi\in\mathbb{F}%
_{s}^{\wedge}}\chi\left(  \operatorname*{Tr}\nolimits_{\mathbb{F}%
_{t}/\mathbb{F}_{s}}\left(  f\left(  u\right)  \right)  \right)  \right) \\
&  =\left(  -1\right)  ^{\rho}t^{2N-\rho}A\left(  s,v\right)  .
\end{align*}
\quad2) Define $f_{a}=af$. Since $a\in\mathbb{F}_{s}$, $f_{a}$ is a hermitian
quadratic form on~$\mathbb{F}_{t}^{2N}$\ and the bilinear form $B_{a}\left(
x,y\right)  =f_{a}\left(  x+y\right)  -f_{a}\left(  x\right)  -f_{a}\left(
y\right)  $ associated to $f_{a}$ satisfies $\operatorname*{Ker}%
B_{a}=\operatorname*{Ker}B$. Hence we can assume that $a=1$ without loss of
generality. Let $v\notin\left(  \operatorname*{Ker}B\right)  ^{\bot}$. The
first part of the Theorem gives $S\left(  f,0\right)  =\left(  -1\right)
^{\rho}t^{2N-\rho}$ hence $S\left(  f,0\right)  \neq0$. Therefore $S\left(
f,v\right)  =0$ if and only if $S\left(  f,v\right)  \overline{S\left(
f,0\right)  }=0$. We have:
\begin{align*}
S\left(  f,v\right)  \overline{S\left(  f,0\right)  }  &  =\sum\limits_{x,y\in
\mathbb{F}_{t}^{2N}}\psi\left(  \left(  f\left(  x\right)  -f\left(  y\right)
+v.x\right)  \right) \\
&  =\sum\limits_{x,y\in\mathbb{F}_{t}^{2N}}\psi\left(  \left(  f\left(
x+y\right)  -f\left(  y\right)  +v.x+v.y\right)  \right) \\
&  =\sum\limits_{x,y\in\mathbb{F}_{t}^{2N}}\psi\left(  \left(  f\left(
x\right)  +B\left(  x,y\right)  +v.x+v.y\right)  \right) \\
&  =\sum\limits_{x\in\mathbb{F}_{t}^{2N}}\psi\left(  \left(  f\left(
x\right)  +v.x\right)  \right)  \sum\limits_{y\in\mathbb{F}_{t}^{2N}}%
\psi\left(  \left(  T\left(  x\right)  +v\right)  .y\right)  .
\end{align*}
Since $v\notin\left(  \operatorname*{Ker}B\right)  ^{\bot}$, the sum $T\left(
x\right)  +v$ is never null and the map $l\left(  y\right)  =\left(  T\left(
x\right)  +v\right)  .y$ is a non trivial linear form on$~\mathbb{F}_{t}^{2N}%
$. We conclude from Lemma~\ref{990925a} that $\sum_{y\in\mathbb{F}_{t}^{2N}%
}\psi\left(  l\left(  y\right)  \right)  =0$, and finally that~$S\left(
f,v\right)  \overline{S\left(  f,0\right)  }=0$.%
\endproof
\medskip

\textbf{Remark:} The constant $A\left(  s,v\right)  $ in Theorem~\ref{cher}
depends wether $\operatorname*{Tr}\nolimits_{\mathbb{F}_{t}/\mathbb{F}_{s}%
}\left(  f\left(  u\right)  \right)  =0$ or not. It has a meaning if we check
that $v=T\left(  u\right)  =T\left(  u^{\prime}\right)  $ and
$\operatorname*{Tr}\nolimits_{\mathbb{F}_{t}/\mathbb{F}_{s}}\left(  f\left(
u\right)  \right)  =0$ imply $\operatorname*{Tr}\nolimits_{\mathbb{F}%
_{t}/\mathbb{F}_{s}}\left(  f\left(  u^{\prime}\right)  \right)  =0$. Let
$v=T\left(  u\right)  =T\left(  u^{\prime}\right)  $. Then $u-u^{\prime}%
:=w\in\operatorname*{Ker}T$ and $B\left(  w,u^{\prime}\right)  =f\left(
u\right)  -f\left(  w\right)  -f\left(  u^{\prime}\right)  $. From $B\left(
w,u^{\prime}\right)  =T\left(  w\right)  .u^{\prime}=0$ and $f\left(
w\right)  =\frac{1}{2}B\left(  w,w\right)  =\frac{1}{2}T\left(  w\right)
.w=0$ it follows that $f\left(  u\right)  =f\left(  u^{\prime}\right)  $,
which gives the desired conclusion.

\section{Number of solutions of some trace equations}

\begin{theorem}
\label{991013d}Let $v\in\mathbb{F}{_{t}}^{2N}$, let $\rho$ be a positive
integer such that $1\leq\rho\leq N$, and $f$ be a quadratic hermitian form of
rank $2\rho$ on$~\mathbb{F}{_{t}}^{2N}$. The number $M$ of solutions of the
equation $\operatorname*{Tr}\nolimits_{\mathbb{F}_{t}/\mathbb{F}_{s}}(f\left(
x\right)  +v.x)=0$ in $\mathbb{F}{_{t}^{2N}}$ is
\[
M=\left\{
\begin{array}
[c]{ll}%
\frac{1}{s}\left(  {t^{2N}+}\left(  -1\right)  ^{\rho}A\left(  s,v\right)
{t^{2N-\rho}}\right)  & \text{if }v\in\left(  \operatorname*{Ker}B\right)
^{\bot}=\operatorname*{Im}T,\\
\frac{{t^{2N}}}{s} & \text{else.}%
\end{array}
\right.
\]
\end{theorem}

\textbf{Proof }: Let us introduce the additive character $\psi^{\prime}\left(
x\right)  =\exp\left(  \frac{i2\pi}{p}\operatorname*{Tr}\nolimits_{\mathbb{F}%
_{s}/\mathbb{F}_{p}}\left(  x\right)  \right)  $ on $\mathbb{F}_{s}$.
Theorem~\ref{991013c} gives
\[
sM=\sum\limits_{c\in\mathbb{F}_{s}}\sum\limits_{x\in\mathbb{F}_{t}^{2N}}%
\psi^{\prime}\left(  c\operatorname*{Tr}\nolimits_{\mathbb{F}_{t}%
/\mathbb{F}_{s}}\left(  f\left(  x\right)  +v.x\right)  \right)  .
\]
Since
\begin{align*}
\psi^{\prime}\left(  c\operatorname*{Tr}\nolimits_{\mathbb{F}_{t}%
/\mathbb{F}_{s}}\left(  f\left(  x\right)  +v.x\right)  \right)   &
=\exp\left(  \frac{i2\pi}{p}\operatorname*{Tr}\nolimits_{\mathbb{F}%
_{s}/\mathbb{F}_{p}}\left(  c\operatorname*{Tr}\nolimits_{\mathbb{F}%
_{t}/\mathbb{F}_{s}}\left(  f\left(  x\right)  +v.x\right)  \right)  \right)
\\
&  =\exp\left(  \frac{i2\pi}{p}\operatorname*{Tr}\nolimits_{\mathbb{F}%
_{t}/\mathbb{F}_{p}}\left(  cf\left(  x\right)  +cv.x\right)  \right) \\
&  =\psi\left(  cf\left(  x\right)  +cv.x\right)  ,
\end{align*}
we deduce
\begin{align*}
sM  &  =\sum\limits_{c\in\mathbb{F}_{s}}\sum\limits_{x\in\mathbb{F}_{t}^{2N}%
}\psi\left(  cf\left(  x\right)  +cv.x\right)  =t^{2N}+\sum\limits_{c\in
\mathbb{F}_{s}^{\ast}}\sum\limits_{x\in\mathbb{F}_{t}^{2N}}\psi\left(
c^{-1}f\left(  cx\right)  +v.\left(  cx\right)  \right) \\
&  =t^{2N}+\sum\limits_{c\in\mathbb{F}_{s}^{\ast}}S\left(  c^{-1}f,v\right)  .
\end{align*}
Now the assertion follows from Theorem~\ref{cher}.%
\endproof

\begin{theorem}
\label{azt}\{\cite{mercier 2001}, Prop. 3\} Let $a$ be an element of
$\mathbb{F}_{s}$, $\rho$ be a positive integer with $1\leq\rho\leq N$, and $f$
be a quadratic hermitian form of rank $2\rho$ on$~\mathbb{F}{_{t}}^{2N}$. The
number $M$ of solutions of the equation $\operatorname*{Tr}%
\nolimits_{\mathbb{F}_{t}/\mathbb{F}_{s}}(f\left(  x\right)  )=a$
in$~\mathbb{F}{_{t}^{2N}}$ is
\[
M=\left\{
\begin{array}
[c]{ll}%
\frac{1}{s}\left(  {t^{2N}-}\left(  -1\right)  ^{\rho}{t^{2N-\rho}}\right)  &
\text{if }a\neq0,\\
\frac{1}{s}\left(  {t^{2N}+}\left(  -1\right)  ^{\rho}\left(  s-1\right)
{t^{2N-\rho}}\right)  & \text{if }a=0.
\end{array}
\right.
\]
\end{theorem}

\textbf{Proof }: We can assume that $f$ is given in the standard form
$f\left(  x\right)  =H\left(  y,y\right)  =y_{1}^{t+1}+...+y_{\rho}^{t+1}$
where $y=\iota\left(  x\right)  \in\mathbb{F}_{t^{2}}^{N}$. If~$\mathbb{F}%
{_{s}^{\symbol{94}}}$ denotes the set of additive characters on$~\mathbb{F}%
_{s}$, then (Theorem~\ref{991013c})
\[
sM=\underset{\psi\in\mathbb{F}{_{s}^{\symbol{94}}}}{\sum}\underset
{x\in\mathbb{F}_{t}^{2N}}{\sum}\psi\left(  \operatorname*{Tr}%
\nolimits_{\mathbb{F}_{t}/\mathbb{F}_{s}}\left(  f\left(  x\right)  \right)
-a\right)  .
\]
Hence
\[
sM=t^{2N}+\underset{\psi\neq\mathbf{1}}{\sum}\overline{\psi\left(  a\right)
}\underset{y\in\mathbb{F}_{t^{2}}^{N}}{\sum}\psi\left(  \operatorname*{Tr}%
\nolimits_{\mathbb{F}_{t}/\mathbb{F}_{s}}\left(  y_{1}^{t+1}+...+y_{\rho
}^{t+1}\right)  \right)  .
\]
We have
\begin{align*}
A_{\psi}  &  =\sum\limits_{y\in\mathbb{F}_{t^{2}}^{N}}\psi\left(
\operatorname*{Tr}\nolimits_{\mathbb{F}_{t}/\mathbb{F}_{s}}\left(  y_{1}%
^{t+1}\right)  \right)  ...\psi\left(  \operatorname*{Tr}\nolimits_{\mathbb{F}%
_{t}/\mathbb{F}_{s}}\left(  y_{\rho}^{t+1}\right)  \right) \\
&  =t^{2\left(  N-\rho\right)  }\left(  \sum\limits_{y\in\mathbb{F}_{t^{2}}%
}\psi\left(  \operatorname*{Tr}\nolimits_{\mathbb{F}_{t}/\mathbb{F}_{s}%
}\left(  y^{t+1}\right)  \right)  \right)  ^{\rho}=t^{2\left(  N-\rho\right)
}B_{\psi}^{\rho}%
\end{align*}
where $B_{\psi}=\sum\nolimits_{y\in\mathbb{F}_{t^{2}}}\psi\left(
\operatorname*{Tr}\nolimits_{\mathbb{F}_{t}/\mathbb{F}_{s}}\left(
y^{t+1}\right)  \right)  $. Since the norm $N:\mathbb{F}_{t^{2}}^{\ast
}\rightarrow\mathbb{F}_{t}^{\ast}$ is surjective and satisfies $\left|
N^{-1}\left(  z\right)  \right|  =t+1$ for all $z\in\mathbb{F}_{t}^{\ast}$
(Lemma~\ref{990921b}), we get
\[
B_{\psi}=1+\left(  t+1\right)  \underset{z\in\mathbb{F}_{t}^{\ast}}{\sum}%
\psi\left(  \operatorname*{Tr}\nolimits_{\mathbb{F}_{t}/\mathbb{F}_{s}}\left(
z\right)  \right)  =-t.
\]
Therefore
\[
sM=t^{2N}+\left(  -1\right)  ^{\rho}t^{2N-\rho}\underset{\psi\neq\mathbf{1}%
}{\sum}\overline{\psi\left(  a\right)  }=t^{2N}+\left(  -1\right)  ^{\rho
}t^{2N-\rho}\left(  -1+\underset{\psi\in\mathbb{F}{_{s}^{\symbol{94}}}}{\sum
}\overline{\psi}\left(  a\right)  \right)
\]
and the usual orthogonality relation establishes the formula.%
\endproof
\medskip

\textbf{Remark} : Theorem~\ref{azt} follows from Theorem~\ref{991013d}
when$~a=0$. A generalization of these two results would be to compute the
number of solutions of $\operatorname*{Tr}\nolimits_{\mathbb{F}_{t}%
/\mathbb{F}_{s}}(f\left(  x\right)  +v.x)=a$ in$~\mathbb{F}{_{t}^{2N}}$.

\section{\label{991014a}The code $\Gamma$}

Remember that $\operatorname*{QH}\left(  \mathbb{F}{_{t}^{2N}}\right)  $
denotes the $\mathbb{F}{_{t}}$-vector space of quadratic hermitian forms
on~$\mathbb{F}{_{t}^{2N}}$. The image of the linear map
\[%
\begin{array}
[c]{cccc}%
\gamma: & \operatorname*{QH}\left(  \mathbb{F}{_{t}^{2N}}\right)
\times\mathbb{F}{_{t}^{2N}} & \rightarrow & \mathbb{F}_{s}^{t^{2N}}\\
& \left(  f,v\right)  & \mapsto & \left(  \operatorname*{Tr}%
\nolimits_{\mathbb{F}_{t}/\mathbb{F}_{s}}(f\left(  x\right)  +v.x)\right)
_{x\in\mathbb{F}{_{t}^{2N}}}%
\end{array}
\]
is a code $\Gamma$ in $\mathbb{F}_{s}^{t^{2N}}$. This code was first
introduced by J.-P. Cherdieu in~\cite{cherdieu} and next~Theorem provides us
whith its parameters. Let us denote by $\operatorname*{w}\left(  c\right)  $
the weight of a non null code-word in a code $C$. If $d\leq\operatorname*{w}%
\left(  c\right)  \leq D$ and if the bounds of these inequalities are reached,
we say that $d$ is the minimal distance of~$C$, and that $r=\frac{D}{d}$ is
the disparity of~$C$.

\begin{theorem}
The weights $\operatorname*{w}\left(  \gamma\left(  f,v\right)  \right)  $ of
the non null code-word $\gamma\left(  f,v\right)  $ of the code$~\Gamma$
satisfy:
\[
t^{2N}-\frac{1}{s}\left(  t^{2N}+t^{2N-1}\right)  \leq\operatorname*{w}\left(
\gamma\left(  f,v\right)  \right)  \leq t^{2N}-\frac{1}{s}\left(
t^{2N}-\left(  s-1\right)  t^{2N-1}\right)
\]
and the bounds of these inequalities are reached. The parameters and the
disparity of$~\Gamma$ are:
\[
\left[  N_{\Gamma},K_{\Gamma},D_{\Gamma}\right]  =\left[  t^{2N},\left(
N^{2}+2N\right)  \log_{s}t,t^{2N}-\frac{1}{s}\left(  t^{2N}+t^{2N-1}\right)
\right]  \text{ and }r\left(  \Gamma\right)  =\frac{\left(  s-1\right)
\left(  t+1\right)  }{st-t-1}.
\]
\end{theorem}

\textbf{Proof }: The length of $\Gamma$ is $N_{\Gamma}=t^{2N}$. It follows
immediately from Theorem~\ref{991013d} that the equation $\operatorname*{Tr}%
\nolimits_{\mathbb{F}_{t}/\mathbb{F}_{s}}(f\left(  x\right)  +v.x)=0$ have
$t^{2N}$ solutions in~$\mathbb{F}{_{t}^{2N}}$ if and only if$~\left(
f,v\right)  =\left(  0,0\right)  $. Consequently the map$~\gamma$ is injective
and
\[
K_{\Gamma}=\dim_{\mathbb{F}_{s}}\Gamma=\dim_{\mathbb{F}_{s}}\left(
\operatorname*{QH}\left(  \mathbb{F}{_{t}^{2N}}\right)  \times\mathbb{F}%
{_{t}^{2N}}\right)  =\left(  N^{2}+2N\right)  \log_{s}t.
\]
We have $\operatorname*{w}\left(  \gamma\left(  f,v\right)  \right)
=t^{2N}-\operatorname*{M}\left(  f,v\right)  $ where the number
$\operatorname*{M}\left(  f,v\right)  $ of solutions of the equation
$\operatorname*{Tr}\nolimits_{\mathbb{F}_{t}/\mathbb{F}_{s}}(f\left(
x\right)  +v.x)=0$ in$~\mathbb{F}{_{t}^{2N}}$ is provided by
Theorem~\ref{991013d}:
\[
\operatorname*{M}\left(  f,v\right)  =\left\{
\begin{array}
[c]{ll}%
\frac{1}{s}\left(  {t^{2N}+}\left(  -1\right)  ^{\rho}A\left(  s,v\right)
{t^{2N-\rho}}\right)  & \text{if }v\in\left(  \operatorname*{Ker}B\right)
^{\bot}=\operatorname*{Im}T,\\
\frac{{t^{2N}}}{s} & \text{else}.
\end{array}
\right.
\]
We consider several cases:\smallskip

\textbf{1.} If $v=0$, then $f\neq0$, and

\quad\quad\textbf{1.1.} If $\rho$ is even, then $2\leq\rho\leq2\left[
\frac{N}{2}\right]  $ and%
\[
\frac{1}{s}\left(  {t^{2N}+}\left(  s-1\right)  {t^{2N-2\left[  \frac{N}%
{2}\right]  }}\right)  \leq\operatorname*{M}\left(  f,0\right)  \leq\frac
{1}{s}\left(  {t^{2N}+}\left(  s-1\right)  {t^{2N-2}}\right)  .\quad
\quad\left(  1\right)
\]
\quad\quad\textbf{1.2.} If $\rho$ is odd, then $1\leq\rho\leq2\left[
\frac{N-1}{2}\right]  +1$ and
\[
\frac{1}{s}\left(  {t^{2N}-}\left(  s-1\right)  {t^{2N-2}}\right)
\leq\operatorname*{M}\left(  f,0\right)  \leq\frac{1}{s}\left(  {t^{2N}%
-}\left(  s-1\right)  {t^{2N-2\left[  \frac{N}{2}\right]  }}\right)
.\quad\quad\left(  2\right)
\]
\textbf{2.} If $v\neq0$,

\quad\quad\textbf{2.1.} If $\rho$ is even and $v\in\left(  \operatorname*{Ker}%
B\right)  ^{\bot}$, then $\rho\neq0$. We get
\[
2\leq\rho\leq2\left[  \frac{N}{2}\right]  \text{ et }\operatorname*{M}\left(
f,v\right)  =\frac{1}{s}\left(  {t^{2N}+}A\left(  s,v\right)  {t^{2N-\rho}%
}\right)  .
\]
We can find a vector $u$ such that $\operatorname*{Tr}\nolimits_{\mathbb{F}%
_{t}/\mathbb{F}_{s}}\left(  f\left(  u\right)  \right)  \neq0$ (indeed
$f\left(  u\right)  =y_{1}^{t+1}+...+y_{\rho}^{t+1}$ in a convenient basis,
and the map $y\mapsto\operatorname*{Tr}\nolimits_{\mathbb{F}_{t}%
/\mathbb{F}_{s}}\left(  y^{t+1}\right)  $ is surjective since
$\operatorname*{Tr}\nolimits_{\mathbb{F}_{t}/\mathbb{F}_{s}}$ are
$\operatorname*{N}_{\mathbb{F}_{t^{2}}/\mathbb{F}_{t}}$ are surjective) thus
there will be~$2$ possible cases:

\quad\quad\quad\quad\textbf{2.1.1.} If $v=T\left(  u\right)  $ with
$\operatorname*{Tr}\nolimits_{\mathbb{F}_{t}/\mathbb{F}_{s}}\left(  f\left(
u\right)  \right)  =0$, then $\operatorname*{M}\left(  f,v\right)  =\frac
{1}{s}\left(  {t^{2N}+}\left(  s-1\right)  {t^{2N-\rho}}\right)  $ and
\[
\frac{1}{s}\left(  {t^{2N}+}\left(  s-1\right)  {t^{2N-2\left[  \frac{N}%
{2}\right]  }}\right)  \leq\operatorname*{M}\left(  f,v\right)  \leq\frac
{1}{s}\left(  {t^{2N}+}\left(  s-1\right)  {t^{2N-2}}\right)  .\quad
\quad\left(  3\right)
\]
\quad\quad\quad\quad\textbf{2.1.2.} If $v=T\left(  u\right)  $ with
$\operatorname*{Tr}\nolimits_{\mathbb{F}_{t}/\mathbb{F}_{s}}\left(  f\left(
u\right)  \right)  \neq0$, then $\operatorname*{M}\left(  f,v\right)
=\frac{1}{s}\left(  {t^{2N}-t^{2N-\rho}}\right)  $ and
\[
\frac{1}{s}\left(  {t^{2N}-t^{2N-2}}\right)  \leq\operatorname*{M}\left(
f,v\right)  \leq\frac{1}{s}\left(  {t^{2N}-t^{2N-2\left[  \frac{N}{2}\right]
}}\right)  .\quad\quad\left(  4\right)
\]
\quad\quad\textbf{2.2.} If $\rho$ is even and $v\notin\left(
\operatorname*{Ker}B\right)  ^{\bot}$, then $\operatorname*{M}\left(
f,v\right)  =\frac{{t^{2N}}}{s}$ belongs to one of the intervals defined by
$\left(  3\right)  $ or$~\left(  4\right)  $.\smallskip

\quad\quad\textbf{2.3.} If $\rho$ is odd and $v\in\left(  \operatorname*{Ker}%
B\right)  ^{\bot}$, then $\operatorname*{M}\left(  f,v\right)  =\frac{1}%
{s}\left(  {t^{2N}-}A\left(  s,v\right)  {t^{2N-\rho}}\right)  .$\smallskip

\quad\quad\quad\quad\textbf{2.3.1.} If $v=T\left(  u\right)  $ with
$\operatorname*{Tr}\nolimits_{\mathbb{F}_{t}/\mathbb{F}_{s}}\left(  f\left(
u\right)  \right)  =0$, then $\operatorname*{M}\left(  f,v\right)  =\frac
{1}{s}\left(  {t^{2N}-}\left(  s-1\right)  {t^{2N-\rho}}\right)  $ and
\[
\frac{1}{s}\left(  {t^{2N}-}\left(  s-1\right)  {t^{2N-1}}\right)
\leq\operatorname*{M}\left(  f,v\right)  \leq\frac{1}{s}\left(  {t^{2N}%
-}\left(  s-1\right)  {t^{2N-2\left[  \frac{N-1}{2}\right]  -1}}\right)
.\quad\quad\left(  5\right)
\]
\quad\quad\quad\quad\textbf{2.3.2.} If $v=T\left(  u\right)  $ with
$\operatorname*{Tr}\nolimits_{\mathbb{F}_{t}/\mathbb{F}_{s}}\left(  f\left(
u\right)  \right)  \neq0$, then $\operatorname*{M}\left(  f,v\right)
=\frac{1}{s}\left(  {t^{2N}+t^{2N-\rho}}\right)  $ and
\[
\frac{1}{s}\left(  {t^{2N}+t^{2N-2\left[  \frac{N-1}{2}\right]  -1}}\right)
\leq\operatorname*{M}\left(  f,v\right)  \leq\frac{1}{s}\left(  {t^{2N}%
+t^{2N-1}}\right)  .\quad\quad\left(  6\right)
\]
\quad\quad\textbf{2.4.} If $\rho$ is odd and $v\notin\left(
\operatorname*{Ker}B\right)  ^{\bot}$, then $\operatorname*{M}\left(
f,v\right)  =\frac{{t^{2N}}}{s}$ belongs to one of the intervals defined by
$\left(  3\right)  $ or $\left(  4\right)  $.\smallskip

It is sufficient to consider the bounds $\left(  1\right)  $ \`{a} $\left(
6\right)  $ to deduce
\[
\frac{1}{s}\left(  t^{2N}-\left(  s-1\right)  t^{2N-1}\right)  \leq
\operatorname*{M}\left(  f,v\right)  \leq\frac{1}{s}\left(  {t^{2N}+t^{2N-1}%
}\right)
\]
for all $\left(  f,v\right)  \in\left(  \operatorname*{QH}\left(
\mathbb{F}{_{t}^{2N}}\right)  \times\mathbb{F}{_{t}^{2N}}\right)  {\backslash
}\left\{  \left(  0,0\right)  \right\}  $. Hence we obtain the bounds of the
weights $\operatorname*{w}\left(  \gamma\left(  f,v\right)  \right)  $.%
\endproof
\medskip

Let $C$ denote a code $\left[  N_{C},K_{C},D_{C}\right]  $. The ratio
$\frac{K_{C}}{N_{C}}$ is called the transmission rate, and the ratio
$\frac{D_{C}}{N_{C}}$ represents the reliability of~$C$. Note that $C$ can
correct $\left[  \frac{D_{C}-1}{2}\right]  $ errors and that
\[
\lambda\left(  C\right)  =\frac{K_{C}}{N_{C}}+\frac{D_{C}}{N_{C}}%
\]
is less than $1+\frac{1}{N_{C}}$ and must be as great as possible.\medskip

The generalized Reed-Muller code $R\left(  r,m\right)  $ of order $r$ on
$\mathbb{F}{_{t}^{m}}$ is described by the code-words $\left(  f\left(
x\right)  \right)  _{x\in\mathbb{F}{_{t}^{m}}}$ where $f$ are polynomials
in~$\mathbb{F}{_{t}}\left[  X_{1},...,X_{m}\right]  $ of total degree less
than$~r$. The dimension of $R\left(  r,m\right)  $ is $C_{m+r}^{r}$ if $r<t$,
and the parameters of $R\left(  2,2N\right)  $ are
\[
\left[  N_{R},K_{R},D_{R}\right]  =\left[  t^{2N},2N^{2}+3N+1,t^{2N}%
-2t^{2N-1}\right]  .
\]
Let us compare $R\left(  2,2N\right)  $ to the code $\Gamma$ with same length
$t^{2N}$ obtained with $s=t$. The code $R\left(  2,2N\right)  $ have a better
transmission rate since
\[
\frac{K_{R}}{N_{R}}-\frac{K_{\Gamma}}{N_{\Gamma}}=\frac{1}{t^{2N}}\left(
N^{2}+N+1\right)
\]
is always positive, but the numbers of corrected errors is better with
$\Gamma$ since
\[
D_{\Gamma}-D_{R}=t^{2N-1}-t^{2N-2}%
\]
is always positive. One can also check that the difference
\[
\lambda\left(  \Gamma\right)  -\lambda\left(  R\right)  =\frac{1}{t^{2N}%
}\left(  t^{2N-1}-t^{2N-2}-N^{2}-N-1\right)
\]
is positive or null as soon as $N\geq2$ or~$t\geq4$. In this sens, $\Gamma$
have better parameters than$~R\left(  2,2N\right)  $.

\section{\label{991101a}The code $C$}

The parameters of the code $\Gamma$ in Section~\ref{991014a} are computed from
Theorem~\ref{991013d}. We can apply the same construction to use
Theorem~\ref{azt}. The image of the linear map
\[%
\begin{array}
[c]{cccc}%
c: & \operatorname*{QH}\left(  \mathbb{F}{_{t}^{2N}}\right)  \times
\mathbb{F}{_{s}} & \rightarrow & \mathbb{F}_{s}^{t^{2N}}\\
& \left(  f,a\right)  & \mapsto & \left(  \operatorname*{Tr}%
\nolimits_{\mathbb{F}_{t}/\mathbb{F}_{s}}(f\left(  x\right)  )-a\right)
_{x\in\mathbb{F}{_{t}^{2N}}}%
\end{array}
\]
is a code $C$ with length $N_{C}=t^{2N}$ on $\mathbb{F}_{s}$. The map $c$ is
one to one. Indeed, if the non null quadratic form $f$ satisfies
$\operatorname*{Tr}\nolimits_{\mathbb{F}_{t}/\mathbb{F}_{s}}(f\left(
x\right)  )=a$ for all $x\in\mathbb{F}{_{t}^{2N}}$, and if $\rho$ denotes the
rank of $f$, then $f\left(  x\right)  =y_{1}^{t+1}+...+y_{\rho}^{t+1}$ where
$y=\iota x$ and $1\leq\rho\leq N$, and the assumpion on $f$ implies
$\operatorname*{Tr}\nolimits_{\mathbb{F}_{t}/\mathbb{F}_{s}}(y_{1}^{t+1})=a$
for all $y_{1}\in\mathbb{F}{_{t^{2}}}$. This is a contradiction of the fact
that the map $\operatorname*{Tr}\nolimits_{\mathbb{F}_{t}/\mathbb{F}_{s}}%
\circ\operatorname*{N}_{\mathbb{F}_{t^{2}}/\mathbb{F}_{t}}:\mathbb{F}_{t^{2}%
}\rightarrow\mathbb{F}_{s}$ is onto.\smallskip

As $c$ is one to one, the dimension of $C$ will be:
\[
K_{C}=\dim_{\mathbb{F}_{s}}\left(  \operatorname*{QH}\left(  \mathbb{F}%
{_{t}^{2N}}\right)  \times\mathbb{F}{_{s}}\right)  {=1+N}^{2}\log_{s}t.
\]

\begin{theorem}
The weights $\operatorname*{w}\left(  c\left(  f,a\right)  \right)  $ of the
non null code-words $c\left(  f,a\right)  $ in $C$ satisfy:
\[
t^{2N}-\frac{1}{s}\left(  t^{2N}+t^{2N-1}\right)  \leq\operatorname*{w}\left(
c\left(  f,a\right)  \right)  \leq t^{2N}%
\]
and the bounds are reached. The parameters and the disparity of $C$ are:
\[
\left[  N_{C},K_{C},D_{C}\right]  =\left[  t^{2N},{1+N}^{2}\log_{s}%
t,t^{2N}-\frac{1}{s}\left(  t^{2N}+t^{2N-1}\right)  \right]  \text{ and
}r\left(  C\right)  =\frac{st}{st-t-1}\text{.}%
\]
\end{theorem}

\textbf{Proof }: It suffices to bound the weights $\operatorname*{w}\left(
c\left(  f,a\right)  \right)  $. We certainly have
\[
\operatorname*{w}\left(  c\left(  f,a\right)  \right)  =t^{2N}%
-\operatorname*{M}\left(  f,a\right)
\]
where $\operatorname*{M}\left(  f,a\right)  $, which denotes the number of
solutions of the equation $\operatorname*{Tr}\nolimits_{\mathbb{F}%
_{t}/\mathbb{F}_{s}}(f\left(  x\right)  )=a$ in $\mathbb{F}{_{t}^{2N}}$, is
given by Theorem~\ref{azt}.\medskip

\textbf{1}. If $a=0$, we know that $\rho\neq0$.\smallskip

\quad\quad\textbf{1}.\textbf{1.} If $\rho$ is even, then $2\leq\rho
\leq2\left[  \frac{N}{2}\right]  $ and
\[
\frac{1}{s}\left(  {t^{2N}+}\left(  s-1\right)  {t^{2N-2\left[  \frac{N}%
{2}\right]  }}\right)  \leq\operatorname*{M}\left(  f,0\right)  \leq\frac
{1}{s}\left(  {t^{2N}+}\left(  s-1\right)  {t^{2N-2}}\right)  .\quad
\quad\left(  1\right)
\]
\quad\quad\textbf{1}.\textbf{2.} If $\rho$ is odd, then $1\leq\rho\leq2\left[
\frac{N-1}{2}\right]  +1$ and
\[
\frac{1}{s}\left(  {t^{2N}-}\left(  s-1\right)  {t^{2N-1}}\right)
\leq\operatorname*{M}\left(  f,0\right)  \leq\frac{1}{s}\left(  {t^{2N}%
-}\left(  s-1\right)  {t^{2N-2\left[  \frac{N-1}{2}\right]  -1}}\right)
.\quad\quad\left(  2\right)
\]
\textbf{2}. If $a\neq0$,

\quad\quad\textbf{2}.\textbf{1.} If $\rho$ is even,
\[
0\leq\operatorname*{M}\left(  f,a\right)  \leq\frac{1}{s}\left(
{t^{2N}-t^{2N-2\left[  \frac{N}{2}\right]  }}\right)  .\quad\quad\left(
3\right)
\]
\quad\quad\textbf{2}.\textbf{2.} If $\rho$ is odd,
\[
\frac{1}{s}\left(  {t^{2N}-t^{2N-2\left[  \frac{N-1}{2}\right]  {-1}}}\right)
\leq\operatorname*{M}\left(  f,a\right)  \leq\frac{1}{s}\left(  {t^{2N}%
+t^{2N-1}}\right)  .\quad\quad\left(  4\right)
\]
The bounds $\left(  1\right)  $ to $\left(  4\right)  $ imply
\[
\forall\left(  f,a\right)  \in\left(  \operatorname*{QH}\left(  \mathbb{F}%
{_{t}^{2N}}\right)  \times\mathbb{F}{_{s}}\right)  {\backslash}\left\{
\left(  0,0\right)  \right\}  {\quad\quad}0\leq\operatorname*{M}\left(
f,a\right)  \leq\frac{1}{s}\left(  {t^{2N}+t^{2N-1}}\right)  .%
\endproof
\]
Let us compare $C$ with $\Gamma$ and $R\left(  2,2N\right)  $. The codes $C$
and $\Gamma$ have same length and same minimal distance, thus will correct the
same amount of errors. Nevertheless the dimension of $\Gamma$ is greater than
those of $C$, hence $\Gamma$ is better at this point of view. But $C$ can be
compared with the Reed-Muller code $R\left(  2,2N\right)  $ when $s=t$. Since
\[
D_{C}-D_{R}=t^{2N-2}\left(  t-1\right)  >0
\]
we find that $C$ can correct more errors than $R\left(  2,2N\right)  $. But
the\ transmission rate is not so good because
\[
\frac{K_{R}}{N_{R}}-\frac{K_{C}}{N_{C}}=\frac{N^{2}+3N}{t^{2N}}>0.
\]
We can check that $\lim\limits_{t\rightarrow+\infty}\left(  \lambda\left(
C\right)  -\lambda\left(  R\right)  \right)  =0$ when $N$ is chosen. In this
sens, $C$ can be compared with $R\left(  2,2N\right)  $ for large values
of$~t$. In the same manner $\lim\limits_{t\rightarrow+\infty}\left(
\frac{K_{R}}{N_{R}}-\frac{K_{C}}{N_{C}}\right)  =0$ and the transmission rates
of $C$ and $R\left(  2,2N\right)  $ can be compared for large values of$~t$.

\section{\label{991013b}Annex: Group characters}

Let $\left(  G,+\right)  $ be a finite abelian group of order$~\left|
G\right|  $. A character $\psi$ on $G$ is a homomorphism from $\left(
G,+\right)  $ to the multiplicative group $\left(  \mathbb{C}^{\ast},.\right)
$ of non null complex numbers. It is easily seen that all $z$ in
$\operatorname*{Im}\psi$ has absolute value~$1$, that $\psi\left(  0\right)
=1$ and $\psi\left(  -x\right)  =\psi\left(  x\right)  ^{-1}=\overline
{\psi\left(  x\right)  }$ for all $x\in G$. The trivial character $\mathbf{1}$
is defined by $\mathbf{1}\left(  x\right)  =1$ for all $x\in G$. The set
$G^{\wedge}$ of all characters defined on$~G$ is a multiplicative group of
order$~\left|  G\right|  $, with the natural law $\left(  \psi\chi\right)
\left(  x\right)  =\psi\left(  x\right)  .\chi\left(  x\right)  $. We have:

\begin{theorem}
\label{991008b}\{\cite{lidl1}, Theorem 5.4\} \textbf{Orthogonality relations
(I).}\newline If $\psi\in G^{\wedge},$%
\[
\sum\limits_{x\in G}\psi\left(  x\right)  =\left\{
\begin{array}
[c]{l}%
0\text{ if }\psi\neq\mathbf{1}\text{,}\\
\left|  G\right|  \text{ else.}%
\end{array}
\right.
\]
If $x\in G$,
\[
\sum\limits_{\psi\in G^{\wedge}}\psi\left(  x\right)  =\left\{
\begin{array}
[c]{l}%
0\text{ if }x\neq0\text{,}\\
\left|  G\right|  \text{ else.}%
\end{array}
\right.
\]
\end{theorem}

Theorem~\ref{991008b} immediately gives us two useful results:

\begin{theorem}
\textbf{Orthogonality relations (II).}
\begin{gather*}
\text{If }\psi,\chi\in G^{\wedge}\text{, then }\sum\limits_{x\in G}\psi\left(
x\right)  \overline{\chi\left(  x\right)  }=\left\{
\begin{array}
[c]{l}%
0\text{ if }\psi\neq\chi\text{,}\\
\left|  G\right|  \text{ else.}%
\end{array}
\right. \\
\text{If }x,y\in G\text{, then }\sum\limits_{\psi\in G^{\wedge}}\psi\left(
x\right)  \overline{\psi\left(  y\right)  }=\left\{
\begin{array}
[c]{l}%
0\text{ if }x\neq y\text{,}\\
\left|  G\right|  \text{ else.}%
\end{array}
\right.
\end{gather*}
\end{theorem}

\textbf{Proof }: We notice that $\psi\left(  x\right)  \overline{\chi\left(
x\right)  }=\left(  \psi\chi^{-1}\right)  \left(  x\right)  $ and $\psi\left(
x\right)  \overline{\psi\left(  y\right)  }=\psi\left(  x-y\right)  $, and we
apply Theorem~\ref{991008b}.%
\endproof

\begin{theorem}
\label{991013c}Let $f:E\rightarrow G$ be a map from a set $E$ to a finite
abelian group$~G$, and let $a\in G$. The number $M$ of solutions of the
equation $f\left(  x\right)  =a$ is
\[
M=\frac{1}{\left|  G\right|  }\sum\limits_{\psi\in G^{\wedge}}\sum
\limits_{x\in E}\psi\left(  f\left(  x\right)  -a\right)  .
\]
\end{theorem}

\textbf{Proof }:
\[
\sum\limits_{\psi\in G^{\wedge}}\sum\limits_{x\in E}\psi\left(  f\left(
x\right)  -a\right)  =\sum\limits_{x\in f^{-1}\left(  a\right)  }%
\sum\limits_{\psi\in G^{\wedge}}\psi\left(  0\right)  +\sum\limits_{x\notin
f^{-1}\left(  a\right)  }\sum\limits_{\psi\in G^{\wedge}}\psi\left(  f\left(
x\right)  -a\right)  =M\left|  G\right|  .%
\endproof
\]
When $G$ is cyclic of order $n$ and generated by $g$, we can check that
$G^{\wedge}$ is cyclic and
\[
G^{\wedge}=\left\{  \psi_{j}\,/\,j\in\left\{  0,1,...,n\right\}  \right\}
\text{ with }\psi_{j}\left(  g^{k}\right)  =\exp\left(  \frac{ijk2\pi}%
{n}\right)  \text{.}%
\]
An additive character on $\mathbb{F}_{t}$ is a character of the additive
group$~\left(  \mathbb{F}_{t},+\right)  $. One can prove that
\[
\psi\left(  x\right)  =\exp\left(  \frac{i2\pi}{p}\operatorname*{Tr}%
\nolimits_{\mathbb{F}_{t}/\mathbb{F}_{p}}\left(  x\right)  \right)
\]
define a non trivial additive character on $\mathbb{F}_{t}$, and that all
others additive characters are given by~ $\psi_{a}\left(  x\right)
=\psi\left(  ax\right)  $ where$~a\in\mathbb{F}_{t}$.

\end{document}